\documentclass[english,fleqn,allpages]{ISTE-Standard}[2005/11/14]

\usepackage{natbib}

\usepackage{amsthm,mathrsfs,verbatim}

\newsavebox{\fminibox} \newlength{\fminilength}

\title{The title to come}
\makeatletter
\let\bbl@nonfrenchlistspacing\relax
\makeatother

\usepackage{amsmath,amssymb,amsthm,graphicx,epic,lscape,pifont,enumerate}
\usepackage{rotating}
\usepackage{multirow}
\usepackage[utf8]{inputenc}
\usepackage{times}
\usepackage{natbib}
\usepackage{xcolor}
\usepackage{pgf}
\usepackage{tikz}
\usepackage{tikz-qtree}
\usepackage{forest}
\usetikzlibrary{shadows,trees,arrows}
\usetikzlibrary{shapes,positioning}
\let\cite=\citet
\bibpunct{[}{]}{,}{n}{}{;}

\numberwithin{equation}{section}
\numberwithin{table}{subsection}
\numberwithin{figure}{subsection}

\usepackage{color}
\definecolor{lightred}{rgb}{1,0.8,0.8}
\definecolor{lightblue}{rgb}{0.7,0.7,1}
\definecolor{lightgreen}{rgb}{0.8,1,0.8}
\definecolor{purple}{rgb}{1,0,1}
\definecolor{grey}{rgb}{0.8,0.8,0.8}

\usepackage{textcomp}

\newcommand{\Nset}{\mathbb{N}}

\newcommand{\set}[1]{\left\{ #1\right\}}

\newcommand{\E}{\mbox{$\mathbb{E}$}}

\newcommand{\PP}{\mbox{$\mathbb{P}$}}
\DeclareMathOperator{\casc}{casc}

\newcommand{\CC}{\mathcal{C}}

\newcommand{\CR}{\mathcal{R}}
\newcommand{\CA}{\mathcal{A}}

\newcommand{\CT}{\mathcal{T}}

\newcommand{\CS}{\mathcal{S}}

\newcommand\1{\leavevmode\hbox{\rm \small1\kern-0.35em\normalsize1}}

\DeclareMathOperator{\pref}{pref}

\newcommand{\ttN}{\mathtt n}
\newcommand{\ttS}{\mathtt s}
\newcommand{\ttE}{\mathtt e}
\newcommand{\ttW}{\mathtt w}

\newcommand{\g}[1]{\mathbb #1}

\newcommand{\rond}[1]{{\mathcal #1}}

\usepackage[draft]{fixme}
\fxusetheme{color}
\FXRegisterAuthor{a}{aa}{A}
\FXRegisterAuthor{pe}{ape}{P}
\FXRegisterAuthor{b}{ab}{B}
\FXRegisterAuthor{f}{af}{F}
\FXRegisterAuthor{n}{an}{N}

\newcommand\bbox{\hfill\rule{2mm}{2mm}}

\newtheorem{assumption}{Assumption}

\begin{document}
\raggedbottom


\tableofcontents

\mainmatter

\ChapterAuthor[Where VLMC and PRW meet]{Variable Length Markov Chains, Persistent Random Walks: a close encounter}{Peggy \Name{C\'enac}, Brigitte \Name{Chauvin}, Fr\'ed\'eric \Name{Paccaut}, Nicolas \Name{Pouyanne}} \label{chap-struct}

\markboth
{Variable Length Markov Chains, persistent random walks: a close encounter}
{Variable Length Markov Chains, persistent random walks: a close encounter}

\section{Introduction}
This is the story of the encounter between two worlds: the world of random walks and the world of Variable Length Markov Chains (VLMC). The meeting point turns around the semi-Markov property of underlying processes.

In a VLMC, unlike fixed order Markov chains, the probability to predict the next symbol depends on a possibly unbounded part of the past, the length of which depends on the past itself. These relevant parts of pasts are called \emph{contexts}. They are stored in a \emph{context tree}. With each context is associated a probability distribution prescribing the conditional probability of the next symbol, given this context.

Variable length Markov chains are now widely used as random models for character strings. They have been introduced in \cite{rissanen/83} to perform data compression. When they have a finite memory, they provide a parsimonious alternative to fixed order Markov chain models, in which the number of parameters to estimate grows exponentially fast with the order; they are also able to capture finer properties of character sequences. When they have infinite memory -- this will be our case of study in this chapter -- they are a tractable way to build non-Markov models and they may be considered as a subclass of ``cha\^ines à liaisons compl\`etes'' (\cite{doeblin/fortet/37}) or ``chains with infinite order'' (\cite{harris/55}).

Variable length Markov chains are used in bioinformatics, linguistics or coding theory to modelize how words grow or to classify words. In bioinformatics, both for protein families and DNA sequences, identifying patterns that have a biological meaning is a crucial issue. Using VLMC as a model enables to quantify the influence of a meaning pattern by giving a transition probability on the following letter of the sequence. In this way, these patterns appear as contexts of a context tree. Notice that their length may be unbounded (\cite{bejerano/yona/01}).

In addition, if the context tree is recognised to be a signature of a family (of proteins say), this gives an efficient statistical method to test whether or not two samples belong to the same family (\cite{busch/etc/09}).

Therefore, estimating a context tree is an issue of interest and many authors (statisticians or not, applied or not) stress the fact that the height of the context tree should not be supposed to be bounded. This is the case in \cite{galves/leonardi/08} where the algorithm \texttt{CONTEXT} is used to estimate an unbounded context tree or in \cite{garivier/leonardi/11}. Furthermore, as explained in \cite{csiszar/talata/06}, the height of the estimated context tree grows with the sample size so that estimating a context tree by assuming \emph{a priori} that its height is bounded is not realistic.

There is a large litterature on constructing efficient estimators of context trees, as well for finite or infinite context trees.
This chapter is not a review of stastistics issues, which would already be relevant for finite memory VLMC. This is a study of the probabilistic properties of infinite memory VLMC as random processes, and more specifically of the main property of interest for such processes: existence and uniqueness of a stationary measure.

As already been said, VLMC are a natural generalisation to infinite memory of Markov chains. It is usual to index a sequence of random variables forming a Markov chain with positive integers and to make the process grow to the right. The main drawback of this habit for infinite memory process is that the sequence of the process is read from left to right whereas the (possibly infinite) sequence giving the past needed to predict the next symbol is read in the context tree from right to left, thus giving rise to confusion and lack of readability. For this reason, in this chapter, the VLMC grows to the left. In this way, both the process sequence and the memory in the context tree are read from left to right.

Classical random walks have \emph{independent} and identically distributed increments.
In the literature, \emph{Persistent} Random Walks, also called  \emph{Goldstein-Kac Random Walks} or  \emph{Correlated Random Walks} refer to random walks having a 
Markov chain of finite order as an increment process.
For such walks, the dynamics of trajectories has a short memory of given length and the random walk itself is not Markovian any more.
What happens whenever the increments depend on a \emph{non bounded} past memory?

Consider a walker on $\g Z$, allowed to increment its trajectory by $-1$ or $1$ at each step of time.
Assume that the probability to keep the current direction $\pm 1$ depends on the time already spent in the said direction -- the distribution of increments acts thus as a reinforcement of the dependency from the past.
More precisely, the process of increments of such a $1$-dimensional random walk is a Markov chain on the set of (right-)infinite words, with variable -- and unbounded -- length memory: a VLMC.
The concerned VLMC is defined in Section~\ref{subsec:PRWdim1}.
It is based on a context tree called \emph{double comb}.
Besides, Section~\ref{subsec:PRWdim2} deals with a $2$-dimensional persistent random walk defined in an analogous manner on $\g Z^2$ by a VLMC based on a context tree called \emph{quadruple comb}.

These random walks which have an unbounded past memory can be seen as a generalization of ``Directionally Reinforced Random Walks (DRRW)'' introduced by
\cite{Mauldin1996}, in the sense that the persistence times are anisotropic ones.
For a $1$-dimensional random walk associated with a double comb, a complete characterization of recurrence and transience, in terms of changing (or not) direction probabilities, is given in Section~\ref{subsec:PRWdim1}.
More precisely, when one of the random times spent in a given direction (the so-called \emph{persistence times}) is an integrable random variable, the recurrence property is equivalent to a classical drift-vanishing.
In all other cases, the walk is transient unless the weight of the tail distributions of both persistent times are equal.
In dimension $2$, sufficient conditions of transience of recurrence are given in Section~\ref{subsec:PRWdim2}.

Actually, because of the very specific form of the underlying driving VLMC, these persistent random walks turn out to be in one-to-one correspondence with so-called \emph{Markov Additif Processes}.
Section~\ref{sec:meetic} is devoted to the close links between persistent random walks, Markov additive processes, semi-Markov chains and VLMC.

In Section~\ref{sec:VLMCdef}, the definition of a general VLMC and a couple of examples are given.
In Section~\ref{sec:PRW}, the persistent random walks (PRW) are defined and known results on their recurrence properties are collected.
In view to the final Section~\ref{sec:meetic} where we show how PRW and VLMC meet through the world of semi-Markov chains, Section~\ref{sec:VLMC} is devoted to results -- together with a heuristic approach -- on the existence and unicity of stationary measures for a VLMC.

\section{Variable Length Markov Chains: definition of the model}
\label{sec:VLMCdef}
Let $\CA$ be a finite set, called the \emph{alphabet}.
In this text, $\CA$ will most often be the standard alphabet $\CA=\set{0,1}$, but also $\CA=\set{d,u}$ (for \emph{down} and \emph{up}) or $\CA=\set{\ttN,\ttE,\ttW,\ttS}$ (for the cardinal directions).
Let
\[
\CR=\set{\alpha\beta\gamma\cdots : \alpha ,\beta ,\gamma ,\cdots\in\CA}
\]
be the set of \emph{right-infinite} words over $\CA$, written by simple concatenation.
A VLMC on~$\CA$, defined below and most often denoted by $\left( U_n\right) _{n\in\Nset}$, is a particular type of $\CR$-valued discrete time Markov chain where:

\vskip -8pt
- the process evolves between time $n$ and time $n+1$ by adding one letter \emph{on the left} of $U_n$;

\vskip -8pt
- the transition probabilities between time $n$ and time $n+1$ depend on a finite - but not bounded - prefix\footnote{In fact, an infinite prefix might be needed in a denumerable number of cases.
} 
of the current word $U_n$.

Giving a formal frame of such a process leads to the following definitions.
For a complete presentation of VLMC, one can also refer to~\cite{cenac/chauvin/paccaut/pouyanne/12}.

\begin{definition}[Context tree]
A \emph{context tree} on $\CA$ is a saturated tree on $\CA$ having an at most countable set of infinite branches.
\end{definition}

As usual, a \emph{tree on $\CA$} is a set $\CT$ of finite words -- namely a subset of $\cup _{n\in\Nset}\CA^n$ -- which contains the empty word $\emptyset$ (the \emph{root} of $\CT$) and which is prefix-stable: for all finite words $u,v$, $uv\in\CT\Longrightarrow u\in\CT$.
The tree $\CT$ is \emph{saturated} whenever any internal node has $\#\left(\CA\right)$ children:
for any finite word $u$ and for any $\alpha\in\CA$, $u\alpha\in\CT\Longrightarrow\left( \forall\beta\in\CA,~u\beta\in\CT\right)$.
A right-infinite word on $\CA$ is an \emph{infinite branch} of $\CT$ when all its finite prefixes belong to $\CT$.

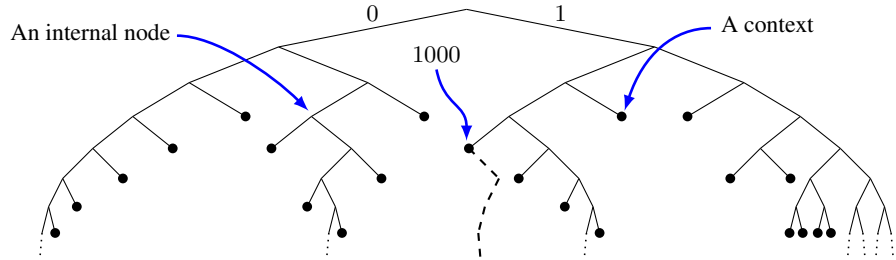
\begin{figure}[h]
\begin{center}
\begin{tikzpicture}[scale=0.5]
\newcommand{\segGeom}[5]{\draw [#5] (0,0) plot[domain=0:#3]({#1+\x*cos(#4)},{#2+\x*sin(#4)});}
\newcommand{\pointCart}[4]{\fill [#4] (0,0) plot [domain=0:360]({#1+#3*cos(\x)},{#2+#3*sin(\x)});}
\newcommand{\pointGeom}[6]{\fill [#6] (0,0) plot [domain=0:360]({#1+#3*cos(#4)+#5*cos(\x)},{#2+#3*sin(#4)+#5*sin(\x)});}
\newcommand{\etA}{1}
\newcommand{\etB}{0.95}%
\newcommand{\etC}{0.9}%
\newcommand{\etD}{0.85}%
\newcommand{\etE}{0.8}%
\newcommand{\etF}{0.75}%
\newcommand{\etG}{0.7}%
\newcommand{\etH}{0.65}%
\newcommand{\penteA}{0.2}
\newcommand{\penteB}{0.4}%
\newcommand{\penteC}{0.6}%
\newcommand{\penteD}{0.8}%
\newcommand{\penteE}{1}%
\newcommand{\penteF}{2}%
\newcommand{\penteG}{4}%
\newcommand{\penteH}{12}%
\newcommand{\noeudB}[2]{\draw [#2] (0,0)--(#1\etA/\penteA,-\etA);}
\newcommand{\noeudC}[3]{\draw [#3] (#1\etA/\penteA,-\etA)--(#1\etA/\penteA#2\etB/\penteB,-\etA-\etB);}
\newcommand{\noeudD}[4]{\draw [#4] (#1\etA/\penteA#2\etB/\penteB,-\etA-\etB)--(#1\etA/\penteA#2\etB/\penteB#3\etC/\penteC,-\etA-\etB-\etC);}
\newcommand{\noeudE}[5]{\draw [#5] (#1\etA/\penteA#2\etB/\penteB#3\etC/\penteC,-\etA-\etB-\etC)--(#1\etA/\penteA#2\etB/\penteB#3\etC/\penteC#4\etD/\penteD,-\etA-\etB-\etC-\etD);}
\newcommand{\noeudF}[6]{\draw [#6] (#1\etA/\penteA#2\etB/\penteB#3\etC/\penteC#4\etD/\penteD,-\etA-\etB-\etC-\etD)--(#1\etA/\penteA#2\etB/\penteB#3\etC/\penteC#4\etD/\penteD#5\etE/\penteE,-\etA-\etB-\etC-\etD-\etE);}
\newcommand{\noeudG}[7]{\draw [#7] (#1\etA/\penteA#2\etB/\penteB#3\etC/\penteC#4\etD/\penteD#5\etE/\penteE,-\etA-\etB-\etC-\etD-\etE)--(#1\etA/\penteA#2\etB/\penteB#3\etC/\penteC#4\etD/\penteD#5\etE/\penteE#6\etF/\penteF,-\etA-\etB-\etC-\etD-\etE-\etF);}
\newcommand{\noeudH}[8]{\draw [#8] (#1\etA/\penteA#2\etB/\penteB#3\etC/\penteC#4\etD/\penteD#5\etE/\penteE#6\etF/\penteF,-\etA-\etB-\etC-\etD-\etE-\etF)--(#1\etA/\penteA#2\etB/\penteB#3\etC/\penteC#4\etD/\penteD#5\etE/\penteE#6\etF/\penteF#7\etG/\penteG,-\etA-\etB-\etC-\etD-\etE-\etF-\etG);}
\newcommand{\noeudI}[9]{\draw [#9] (#1\etA/\penteA#2\etB/\penteB#3\etC/\penteC#4\etD/\penteD#5\etE/\penteE#6\etF/\penteF#7\etG/\penteG,-\etA-\etB-\etC-\etD-\etE-\etF-\etG)--(#1\etA/\penteA#2\etB/\penteB#3\etC/\penteC#4\etD/\penteD#5\etE/\penteE#6\etF/\penteF#7\etG/\penteG#8\etH/\penteH,-\etA-\etB-\etC-\etD-\etE-\etF-\etG-\etH);}
\newcommand{\rayon}{0.13}
\newcommand{\feuilleB}[2]{\draw [#2] (0,0)--(#1\etA/\penteA,-\etA);\fill (#1\etA/\penteA,-\etA) circle(\rayon);}%
\newcommand{\feuilleC}[3]{\draw [#3] (#1\etA/\penteA,-\etA)--(#1\etA/\penteA#2\etB/\penteB,-\etA-\etB);\fill (#1\etA/\penteA#2\etB/\penteB,-\etA-\etB) circle(\rayon);}
\newcommand{\feuilleD}[4]{\draw [#4] (#1\etA/\penteA#2\etB/\penteB,-\etA-\etB)--(#1\etA/\penteA#2\etB/\penteB#3\etC/\penteC,-\etA-\etB-\etC);\fill (#1\etA/\penteA#2\etB/\penteB#3\etC/\penteC,-\etA-\etB-\etC) circle(\rayon);}
\newcommand{\feuilleE}[5]{\draw [#5] (#1\etA/\penteA#2\etB/\penteB#3\etC/\penteC,-\etA-\etB-\etC)--(#1\etA/\penteA#2\etB/\penteB#3\etC/\penteC#4\etD/\penteD,-\etA-\etB-\etC-\etD);\fill (#1\etA/\penteA#2\etB/\penteB#3\etC/\penteC#4\etD/\penteD,-\etA-\etB-\etC-\etD) circle(\rayon);}
\newcommand{\feuilleF}[6]{\draw [#6] (#1\etA/\penteA#2\etB/\penteB#3\etC/\penteC#4\etD/\penteD,-\etA-\etB-\etC-\etD)--(#1\etA/\penteA#2\etB/\penteB#3\etC/\penteC#4\etD/\penteD#5\etE/\penteE,-\etA-\etB-\etC-\etD-\etE);\fill (#1\etA/\penteA#2\etB/\penteB#3\etC/\penteC#4\etD/\penteD#5\etE/\penteE,-\etA-\etB-\etC-\etD-\etE) circle(\rayon);}
\newcommand{\feuilleG}[7]{\draw [#7] (#1\etA/\penteA#2\etB/\penteB#3\etC/\penteC#4\etD/\penteD#5\etE/\penteE,-\etA-\etB-\etC-\etD-\etE)--(#1\etA/\penteA#2\etB/\penteB#3\etC/\penteC#4\etD/\penteD#5\etE/\penteE#6\etF/\penteF,-\etA-\etB-\etC-\etD-\etE-\etF);\fill (#1\etA/\penteA#2\etB/\penteB#3\etC/\penteC#4\etD/\penteD#5\etE/\penteE#6\etF/\penteF,-\etA-\etB-\etC-\etD-\etE-\etF) circle(\rayon);}
\newcommand{\feuilleH}[8]{\draw [#8] (#1\etA/\penteA#2\etB/\penteB#3\etC/\penteC#4\etD/\penteD#5\etE/\penteE#6\etF/\penteF,-\etA-\etB-\etC-\etD-\etE-\etF)--(#1\etA/\penteA#2\etB/\penteB#3\etC/\penteC#4\etD/\penteD#5\etE/\penteE#6\etF/\penteF#7\etG/\penteG,-\etA-\etB-\etC-\etD-\etE-\etF-\etG);\fill (#1\etA/\penteA#2\etB/\penteB#3\etC/\penteC#4\etD/\penteD#5\etE/\penteE#6\etF/\penteF#7\etG/\penteG,-\etA-\etB-\etC-\etD-\etE-\etF-\etG) circle(\rayon);}
\newcommand{\feuilleI}[9]{\draw [#9] (#1\etA/\penteA#2\etB/\penteB#3\etC/\penteC#4\etD/\penteD#5\etE/\penteE#6\etF/\penteF#7\etG/\penteG,-\etA-\etB-\etC-\etD-\etE-\etF-\etG)--(#1\etA/\penteA#2\etB/\penteB#3\etC/\penteC#4\etD/\penteD#5\etE/\penteE#6\etF/\penteF#7\etG/\penteG#8\etH/\penteH,-\etA-\etB-\etC-\etD-\etE-\etF-\etG-\etH);}
\noeudB-{}
\noeudB+{}
\noeudC--{}
\noeudC-+{}
\noeudC+-{}
\noeudC++{}
\noeudD---{}
\feuilleD--+{}
\noeudD-+-{}
\feuilleD-++{}
\noeudD+--{}
\feuilleD+-+{}
\feuilleD++-{}
\noeudD+++{}
\noeudE----{}
\feuilleE---+{}
\feuilleE-+--{}
\noeudE-+-+{}
\feuilleE+---{}\noeudF+---+{dashed,line width=0.8pt}\noeudG+---+-{dashed,line width=0.8pt}\noeudH+---+--{dashed,line width=0.8pt}\noeudI+---+--+{dashed,line width=0.8pt}

\noeudE+--+{}
\noeudE+++-{}
\noeudE++++{}
\noeudF-----{}
\feuilleF----+{}
\noeudF-+-+-{}
\feuilleF-+-++{}
\feuilleF+--+-{}
\noeudF+--++{}
\feuilleF+++--{}
\feuilleF+++-+{}
\noeudF++++-{}
\noeudF+++++{}
\noeudG------{}
\feuilleG-----+{}
\feuilleG-+-+--{}
\noeudG-+-+-+{}
\feuilleG+--++-{}
\noeudG+--+++{}
\noeudG++++--{}
\noeudG++++-+{}
\noeudG+++++-{}
\noeudG++++++{}
\noeudH-------{}
\feuilleH------+{}
\noeudH-+-+-+-{}
\feuilleH-+-+-++{}
\noeudH+--+++-{}
\feuilleH+--++++{}
\feuilleH++++---{}
\feuilleH++++--+{}
\feuilleH++++-+-{}
\feuilleH++++-++{}
\noeudH+++++--{}
\noeudH+++++-+{}
\noeudH++++++-{}
\noeudH+++++++{}
\noeudI--------{dotted, line width = 0.6pt}
\noeudI-+-+-+--{dotted, line width = 0.6pt}
\noeudI+--+++--{dotted, line width = 0.6pt}
\noeudI+++++---{dotted, line width = 0.6pt}
\noeudI+++++-++{dotted, line width = 0.6pt}
\noeudI++++++--{dotted, line width = 0.6pt}
\noeudI++++++++{dotted, line width = 0.6pt}
\draw (-10,-0.6) node{{\small An internal node}};
\draw  [->,>=latex,blue,line width=1pt] (-7.7,-0.7) ..controls +(1,0) and +(-1,1).. (-4.2,-2.7);
\draw (8,-0.4) node{{\small A context}};
\draw  [->,>=latex,blue,line width=1pt] (6.5,-0.5) ..controls +(-1,0) and +(0.5,1.5).. (4.2,-2.7);
\draw (-0.8,-1.2) node{{\small $1000$}};
\draw  [->,>=latex,blue,line width=1pt] (-0.8,-1.6) ..controls +(0.2,-1) and +(0,1).. (0,-3.5);
\draw (-2.5,-0.1) node{{\small $0$}};\draw (2.5,-0.1) node{{\small $1$}};
\end{tikzpicture}
\end{center}
\caption{\label{fig:exampleArbre}A context tree on the alphabet $\CA=\set{0,1}$.
The dotted lines are possibly the beginning of infinite branches.
Any word that writes $1000\cdots$, as the one drawn dashed, admits $1000$ as a $\pref$.}
\end{figure}

\vskip -8pt
A context tree is therefore made of \emph{internal nodes} ($u\in\CT$ is internal when $\exists\alpha\in\CA$, $u\alpha\in\CT$) and of \emph{leaves} ($u\in\CT$ is a leaf when it has no child: $\forall\alpha\in\CA$, $u\alpha\notin\CT$).
Following the vocabulary introduced by Rissanen, a \emph{context} of the tree is a leaf or an infinite branch.
A finite or right-infinite word on $\CA$ is an \emph{external node} when it is neither internal nor a context.
See below Figure~\ref{fig:exampleArbre} that illustrates these definitions, as well as the $\pref$ function defined hereunder.

\begin{definition}[$\pref$ function]
Let $\CT$ be a context tree.
If $w$ is any external node or any context, the symbol $\pref w$ denotes the longest (finite or infinite) prefix of $w$ that belongs to $\CT$.
\end{definition}

In other words, $\pref w$ is the only context $c$ for which $w=c\cdots$
For a more visual presentation, hang $w$ by its head (its left-most letter) and insert it into the tree;
the only context through which the word goes out of the tree is its $\pref$.

With these definitions, it is now possible to define a VLMC.

\begin{definition}[VLMC]
\label{defVLMC}
Let $\CT$ be a context tree.
For every context $c$ of $\CT$, let $q_c$ be a probability measure on $\CA$.
The \emph{variable length Markov chain (VLMC)} defined by $\CT$ and by the $\left( q_c\right) _c$ is the $\CR$-valued discrete-time Markov chain $\left( U_n\right) _{n\in\Nset}$ defined by the following transition probabilities:
$\forall n\in\Nset$, $\forall\alpha\in\CA$,
\begin{equation}
\label{probaTransition}
\PP\left( U_{n+1}=\alpha U_n | U_n\right)
=q_{\pref\left(U_n\right)}\left(\alpha\right).
\end{equation}
\end{definition}

To get a realisation of a VLMC as a process on $\CR$, take a (random) right infinite word
\[
U_0=X_{0}X_{-1}X_{-2}X_{-3}\cdots
\]
At each step of time $n\geq 0$, one gets $U_{n+1}$ by adding a random letter $X_{n+1}$ on the left of $U_n$:
\[
\setlength{\jot}{50pt}
\begin{array}{rl}
U_{n+1}&=X_{n+1}U_n\\[5pt]
&=X_{n+1}X_{n}\cdots X_{1}X_{0}X_{-1}X_{-2}\cdots
\end{array}
\]
under the conditional distribution~\ref{probaTransition}.

\begin{remark}
\emph{Probabilizing} a context tree consists, as in Definition~\ref{defVLMC}, in endowing it with a family of probability measures on the alphabet, indexed by the set of contexts. This vocabulary is used below.
\end{remark}

\begin{remark}
Assume that the context tree is finite and denote its height by $h$;
in this condition, the VLMC is just a Markov chain of order $h$ on $\CA$.
On the contrary, when the context tree is infinite, and this is mainly our case of interest, the VLMC is generally \emph{not} a Markov process on~$\CA$.
\end{remark}

\begin{example}
Take $\CA=\set{\ttN,\ttE,\ttW,\ttS}$ as an (ordered) alphabet, so that the daughters of an internal node are represented as at the left side of Figure~\ref{fig:news}.
Making the transition probabilities $\PP\left( U_{n+1}=\alpha U_n | U_n\right)$ depend only on the length of the largest prefix of the form $\ttN^k$ ($k\geq 0$) of $U_n$ amounts to taking a comb as a context tree, as drawn at the right side of Figure~\ref{fig:news}.
Its finite contexts are the $\ttN^k\alpha$ where $k\geq 0$ and $\alpha\in\CA\setminus\set\ttN$.
\begin{figure}[h]
\begin{center}
\begin{minipage}{150pt}
\begin{tikzpicture}[scale=1.6]
\newcommand{\haut}{0.8}
\newcommand{\ray}{0.05}
\draw (0,0)--(-1/3,-\haut/3);\draw (-2/3,-2*\haut/3)--(-1,-\haut);\fill (-1,-\haut) circle(\ray);\draw (-0.5,-\haut/2) node{$\ttN$};
\draw (0,0)--(-0.3/3,-\haut/3);\draw (-2*0.3/3,-2*\haut/3)--(-0.3,-\haut);\fill (-0.3,-\haut) circle(\ray);\draw (-0.5*0.3,-\haut/2) node{$\ttE$};
\draw (0,0)--(0.3/3,-\haut/3);\draw (2*0.3/3,-2*\haut/3)--(0.3,-\haut);\fill (0.3,-\haut) circle(\ray);\draw (0.5*0.3,-\haut/2) node{$\ttW$};
\draw (0,0)--(1/3,-\haut/3);\draw (2/3,-2*\haut/3)--(1,-\haut);\fill (1,-\haut) circle(\ray);\draw (0.5,-\haut/2) node{$\ttS$};
\end{tikzpicture}
\end{minipage}
\begin{minipage}{100pt}
\includegraphics{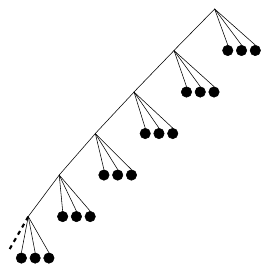}
\end{minipage}
\end{center}
\caption{\label{fig:news}
On the left:
how one can represent trees on $\CA=\set{\ttN,\ttE,\ttW,\ttS}$.
On the right, the so-called \emph{left comb} on $\CA=\set{\ttN,\ttE,\ttW,\ttS}$.}
\end{figure}

\end{example}

\begin{example}
\label{ex:comb}
Take again $\CA=\set{\ttN,\ttE,\ttW,\ttS}$ as an alphabet.
Making the transition probabilities $\PP\left( U_{n+1}=\alpha U_n | U_n\right)$ depend only on the length of the largest prefix of the form $\alpha^k$ ($k\geq 1$) of $U_n$ where $\alpha$ is \emph{any} letter amounts to taking a \emph{quadruple comb} as a context tree, as drawn at the right side of Figure~\ref{fig:bcomb}.
In the same vein, if one takes $\CA=\set{u,d}$, the \emph{double comb} is the context tree drawn at the left side of Figure~\ref{fig:bcomb}.
In the corresponding VLMC, the transitions depend only on the length of the last current run $u^k$ or $d^k$, $k\geq 1$.
The double comb and the quadruple comb are used below to define persistent random walks.
\begin{figure}[h]
\begin{center}
\begin{minipage}{100pt}
\begin{tikzpicture}[scale=0.45]
\newcommand{\segGeom}[5]{\draw [#5] (0,0) plot[domain=0:#3]({#1+\x*cos(#4)},{#2+\x*sin(#4)});}
\newcommand{\pointCart}[4]{\fill [#4] (0,0) plot [domain=0:360]({#1+#3*cos(\x)},{#2+#3*sin(\x)});}
\newcommand{\pointGeom}[6]{\fill [#6] (0,0) plot [domain=0:360]({#1+#3*cos(#4)+#5*cos(\x)},{#2+#3*sin(#4)+#5*sin(\x)});}
\newcommand{\gdePenteG}{0.8}
\newcommand{\gdePenteD}{-\gdePenteG}
\newcommand{\noeud}[2]{-#1*#2,-#1}
\newcommand{\absNoeud}[2]{-#1*#2}
\newcommand{\titePenteGD}{0.3}
\newcommand{\titePenteDG}{-\titePenteGD}
\newcommand{\rayon}{0.13}
\newcommand{\feuille}[3]{
\draw (\noeud{#2}{#1})--++(-#3,-1);
\fill (\absNoeud{#2}{#1}-#3,-#2-1) circle(\rayon);
}
\draw (0,0)--(\noeud{1}{\gdePenteG});
\draw (0,0)--(\noeud{1}{\gdePenteD});
\draw (\noeud{1}{\gdePenteG})--(\noeud{2}{\gdePenteG});
\feuille{\gdePenteG}{1}{\titePenteDG}
\draw (\noeud{1}{\gdePenteD})--(\noeud{2}{\gdePenteD});
\feuille{\gdePenteD}{1}{\titePenteGD}
\draw (\noeud{2}{\gdePenteG})--(\noeud{3}{\gdePenteG});
\feuille{\gdePenteG}{2}{\titePenteDG}
\draw (\noeud{2}{\gdePenteD})--(\noeud{3}{\gdePenteD});
\feuille{\gdePenteD}{2}{\titePenteGD}
\draw (\noeud{3}{\gdePenteG})--(\noeud{3.4}{\gdePenteG});
\draw [dashed] (\noeud{3.4}{\gdePenteG})--(\noeud{4.2}{\gdePenteG});
\feuille{\gdePenteG}{3}{\titePenteDG}
\draw (\noeud{3}{\gdePenteD})--(\noeud{3.4}{\gdePenteD});
\draw [dashed] (\noeud{3.4}{\gdePenteD})--(\noeud{4.2}{\gdePenteD});
\feuille{\gdePenteD}{3}{\titePenteGD}
\end{tikzpicture}
\end{minipage}
\hskip 10pt
\begin{minipage}{200pt}
\begin{tikzpicture}[scale=0.5]
\newcommand{\gdePenteN}{1.8}
\newcommand{\gdePenteE}{0.6}
\newcommand{\gdePenteW}{-\gdePenteE}
\newcommand{\gdePenteS}{-\gdePenteN}
\newcommand{\noeud}[2]{-#1*#2,-#1}
\newcommand{\absNoeud}[2]{-#1*#2}
\newcommand{\titePenteNE}{1}
\newcommand{\titePenteNW}{0.6}
\newcommand{\titePenteNS}{0.2}
\newcommand{\titePenteEN}{0.9}
\newcommand{\titePenteEW}{0.1}
\newcommand{\titePenteES}{-0.3}
\newcommand{\titePenteWN}{-\titePenteES}
\newcommand{\titePenteWE}{-\titePenteEW}
\newcommand{\titePenteWS}{-\titePenteEN}
\newcommand{\titePenteSN}{-\titePenteNS}
\newcommand{\titePenteSE}{-\titePenteNW}
\newcommand{\titePenteSW}{-\titePenteNE}
\newcommand{\rayon}{0.1}
\newcommand{\feuille}[3]{
\draw (\noeud{#2}{#1})--++(-#3,-1);
\fill (\absNoeud{#2}{#1}-#3,-#2-1) circle(\rayon);
}
\draw (0,0)--(\noeud{3.4}{\gdePenteN});
\draw (0,0)--(\noeud{3.4}{\gdePenteE});
\draw (0,0)--(\noeud{3.4}{\gdePenteW});
\draw (0,0)--(\noeud{3.4}{\gdePenteS});
\draw [dashed] (\noeud{3.4}{\gdePenteN})--(\noeud{4.1}{\gdePenteN});
\draw [dashed] (\noeud{3.4}{\gdePenteE})--(\noeud{4.3}{\gdePenteE});
\draw [dashed] (\noeud{3.4}{\gdePenteW})--(\noeud{4.3}{\gdePenteW});
\draw [dashed] (\noeud{3.4}{\gdePenteS})--(\noeud{4.1}{\gdePenteS});
\feuille{\gdePenteN}{1}{\titePenteNE}
\feuille{\gdePenteN}{1}{\titePenteNW}
\feuille{\gdePenteN}{1}{\titePenteNS}
\feuille{\gdePenteE}{1}{\titePenteEN}
\feuille{\gdePenteE}{1}{\titePenteEW}
\feuille{\gdePenteE}{1}{\titePenteES}
\feuille{\gdePenteW}{1}{\titePenteWN}
\feuille{\gdePenteW}{1}{\titePenteWE}
\feuille{\gdePenteW}{1}{\titePenteWS}
\feuille{\gdePenteS}{1}{\titePenteSN}
\feuille{\gdePenteS}{1}{\titePenteSE}
\feuille{\gdePenteS}{1}{\titePenteSW}
\feuille{\gdePenteN}{2}{\titePenteNE}
\feuille{\gdePenteN}{2}{\titePenteNW}
\feuille{\gdePenteN}{2}{\titePenteNS}
\feuille{\gdePenteE}{2}{\titePenteEN}
\feuille{\gdePenteE}{2}{\titePenteEW}
\feuille{\gdePenteE}{2}{\titePenteES}
\feuille{\gdePenteW}{2}{\titePenteWN}
\feuille{\gdePenteW}{2}{\titePenteWE}
\feuille{\gdePenteW}{2}{\titePenteWS}
\feuille{\gdePenteS}{2}{\titePenteSN}
\feuille{\gdePenteS}{2}{\titePenteSE}
\feuille{\gdePenteS}{2}{\titePenteSW}
\feuille{\gdePenteN}{3}{\titePenteNE}
\feuille{\gdePenteN}{3}{\titePenteNW}
\feuille{\gdePenteN}{3}{\titePenteNS}
\feuille{\gdePenteE}{3}{\titePenteEN}
\feuille{\gdePenteE}{3}{\titePenteEW}
\feuille{\gdePenteE}{3}{\titePenteES}
\feuille{\gdePenteW}{3}{\titePenteWN}
\feuille{\gdePenteW}{3}{\titePenteWE}
\feuille{\gdePenteW}{3}{\titePenteWS}
\feuille{\gdePenteS}{3}{\titePenteSN}
\feuille{\gdePenteS}{3}{\titePenteSE}
\feuille{\gdePenteS}{3}{\titePenteSW}
\end{tikzpicture}
\end{minipage}
\end{center}
\caption{\label{fig:bcomb}
The double comb and the quadruple comb.}
\end{figure}
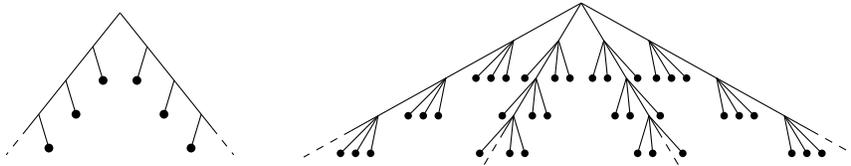
\end{example}

\begin{example}
\label{ex:4VLMC}
Take $\CA=\set{0,1}$ (naturally ordered for the drawings).
The left comb of right combs, drawn at the left side of Figure~\ref{fig:combsComb}, is the context tree of a VLMC that makes its transition probabilities depend on the largest prefix of~$U_n$ of the form~$0^p1^q$.
If one has to take into consideration the largest prefix
of the form $0^p1^q$ or $1^p0^q$, one has to use the double comb of opposite combs, as drawn at the right side of Figure~\ref{fig:combsComb}.
\begin{figure}[h]
\begin{center}
\hbox{
\includegraphics[width=120pt]{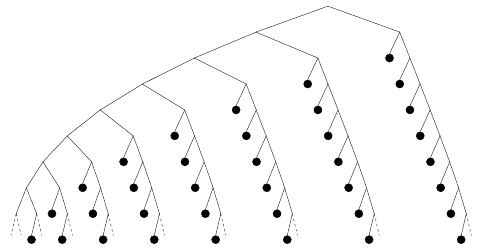}
\hskip 10pt
\includegraphics[width=200pt]{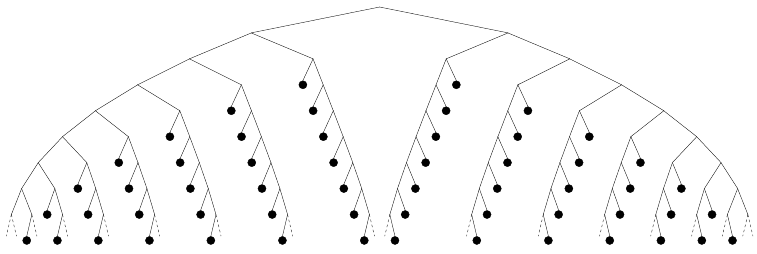}
}
\end{center}
\caption{\label{fig:combsComb}
Context trees on $\CA=\set{0,1}$:
the left comb of right combs (on the left) and a double comb of opposite combs (on the right).}
\end{figure}

\end{example}

\begin{definition}[Non-nullness]
	\label{nn}
A VLMC is called \emph{non-null} when no transition probability vanish, \emph{i.e.} when $q_c(\alpha )>0$ for every context $c$ and for every $\alpha\in\CA$.
\end{definition}

Non-nullness appears below as an irreducibility-like assumption made on the driving VLMC of persistent random walks and for existence and unicity of an invariant probability measure for a general VLMC as well.

\section{Definition and behaviour of Persistent Random Walks}
\label{sec:PRW}

In this section, the so called \emph{Persistent Random Walks (PRW)} are defined.
A PRW is a random walk driven by some VLMC.
In dimension $1$ and $2$, results on transience and recurrence of PRW are given.
These results are detailed and proven in
\cite{CDLO,cenac/chauvin/herrmann/vallois/13}
in dimension one and in
\cite{cenac:hal-01658494}
in dimension two.

\subsection{Persistent Random Walks in dimension one}
\label{subsec:PRWdim1}

In this section, we deal with $1$-dimensional Persistent Random Walks (PRW).
Notice that, contrary to the classical random walk, a PRW is generally not Markovian.
Let $\CA:=\{d,u\}=\{-1,1\}$ ($d$ for down and $u$ for up) and  consider the \emph{double comb} on this alphabet as a context tree, probabilize it and denote by $(U_n)_n$ a realisation of the associated VLMC.
The $n\textsuperscript{th}$ increment $X_n$ of the PRW is given as the first letter of~$U_n$:
define the persistent random walk $S=(S_n)_{n\geq 0}$ by $S_0=0$ and, for $n\geq 1$,
\begin{equation}
\label{rw-S}
S_n:=\sum_{\ell=1}^n X_{\ell},
\end{equation}
so that for any $n\geq 1$, $m\geq 0$,
\begin{eqnarray*}
	\PP\left(S_{m+1}=S_m+1|U_m=d^n u\ldots \right)&=&q_{d^nu}(u) \\
	\PP\left(S_{m+1}=S_m-1|U_m=u^n d\ldots \right)&=&q_{u^nd}(d).
\end{eqnarray*}
Furthermore, for sake of simplicity and without loss of generality, 
we condition the walk to start a.s. from  $\{X_{-1}=u, X_0=d\}$
-- this amounts to changing the origin of time.
In this model, a walker on a line keeps the same direction with a probability which depends on the discrete time already spent in the direction the walker is currently moving.  
See Figure~\ref{marche}.
This model can be seen as a generalisation of Directionally Reinforced Random Walks (DRRWs) introduced in \cite{Mauldin1996}.
 
Taking different probabilized context trees would lead to different probabilistic impacts on the asymptotic behaviour of resulting PRWs. Moreover, the characterization of the recurrent \emph{versus} transient behaviour is difficult in general. We state here exhaustive recurrence criteria for PRWs defined from a double comb. 

In order to avoid trivial cases, we assume that $S$ cannot be frozen in one of the two directions with a positive probability. Therefore, we make the following assumption.
\begin{assumption}[finiteness of the length of runs]
\label{ass:a1}
For any $\alpha,\beta\in\{u,d\}$, $\alpha\neq \beta$, 
\begin{equation}
\label{a1}
\lim _{n\to +\infty}\left(\prod_{k=1}^{n}q_{\alpha^k\beta}(\alpha)\right)=0.
\end{equation}
\end{assumption}

Let $\tau^u_{n}$ and  $\tau^d_{n}$ be respectively  the length of the $n\textsuperscript{th}$ rise and of the $n\textsuperscript{th}$  descent.

\begin{figure}[h]
\definecolor{qqqqff}{rgb}{0,0,1}
\definecolor{cqcqcq}{rgb}{0.75,0.75,0.75}
\begin{center}
\begin{tikzpicture}[scale=0.6,line cap=round,line join=round,>=latex,x=1.0cm,y=1.0cm]
\begin{scriptsize}
\draw [color=cqcqcq,dash pattern=on 3pt off 3pt, xstep=2.0cm,ystep=2.0cm] (-4.87,-3.74) grid (13.61,6.62);
\draw[->,color=black] (-4.87,0) -- (13.61,0);
\foreach \x in {-4,-2,2,4,6,8,10,12}
\draw[shift={(\x,0)},color=black] (0pt,2pt) -- (0pt,-2pt);
\draw[->,color=black] (0,-3.74) -- (0,6.62);
\foreach \y in {-2,2,4,6}
\draw[shift={(0,\y)},color=black] (2pt,0pt) -- (-2pt,0pt);
\clip(-4.87,-3.74) rectangle (13.61,6.62);
\draw [dotted] (-4,0)-- (-3,-1);
\draw [dotted] (-3,-1)-- (-2,-2);
\draw [dotted] (-2,-2)-- (-1,-1);
\draw [dotted] (-1,-1)-- (0,0);
\draw [dotted] (-4,0)-- (-4.26,0.3);
\draw (-3.53,-0.07) node[anchor=north west] {d};
\draw (-2.64,-0.96) node[anchor=north west] {d};
\draw (-1.79,-1.05) node[anchor=north west] {u};
\draw (-0.96,-0.22) node[anchor=north west] {u};
\draw (0,0)-- (1,1);
\draw (1,1)-- (2,2);
\draw (2,2)-- (3,3);
\draw (3,3)-- (4,4);
\draw (4,4)-- (5,3);
\draw (0.15,0.95) node[anchor=north west] {u};
\draw (1.12,1.94) node[anchor=north west] {...};
\draw (2.11,2.93) node[anchor=north west] {...};
\draw (2.99,3.84) node[anchor=north west] {u};
\draw (4.64,3.84) node[anchor=north west] {d};
\draw (3.68,0) node[anchor=north west] {$B_0$};
\draw [->,line width=1pt] (3,4) -- (9,4);
\draw [->,line width=1pt] (4,3) -- (4,6);
\draw [dash pattern=on 6pt off 6pt] (4,3)-- (4,0);
\draw (5,3)-- (6,2);
\draw (6,2)-- (7,3);
\draw [dash pattern=on 6pt off 6pt] (6,2)-- (6,0);
\draw (5.62,0) node[anchor=north west] {$B_1$};
\draw (5.49,2.9) node[anchor=north west] {d};
\draw (6.1,2.85) node[anchor=north west] {u};
\draw (7,3)-- (8,2);
\draw (7.49,2.92) node[anchor=north west] {d};
\draw (6.61,0) node[anchor=north west] {$B_2$};
\draw [dash pattern=on 6pt off 6pt] (7,3)-- (7,0);
\draw (8,2)-- (9,1);
\draw (9,1)-- (10,2);
\draw (10,2)-- (11,3);
\draw (11,3)-- (12,4);
\draw (12,4)-- (12.4,3.6);
\draw (8.44,1.95) node[anchor=north west] {d};
\draw (9.15,1.93) node[anchor=north west] {u};
\draw (10.11,2.9) node[anchor=north west] {u};
\draw (11.1,3.85) node[anchor=north west] {u};
\draw (8.63,0) node[anchor=north west] {$B_3$};
\draw (11.63,0) node[anchor=north west] {$B_4$};
\draw [dash pattern=on 6pt off 6pt] (9,1)-- (9,0);
\draw [dash pattern=on 6pt off 6pt] (12,4)-- (12,0);
\draw [->,line width=1pt] (4,-1) -- (6,-1);
\draw [->,line width=1pt] (6,-1) -- (4,-1);
\draw [->,line width=1pt] (6,-1) -- (7,-1);
\draw [->,line width=1pt] (7,-1) -- (6,-1);
\draw [->,line width=1pt] (7,-1) -- (9,-1);
\draw [->,line width=1pt] (9,-1) -- (7,-1);
\draw [->,line width=1pt] (9,-1) -- (12,-1);
\draw [->,line width=1pt] (12,-1) -- (9,-1);
\draw (4.51,-0.91) node[anchor=north west] {$\tau_1^d$};
\draw (6.03,-1.02) node[anchor=north west] {$\tau_1^u$};
\draw (7.4,-0.90) node[anchor=north west] {$\tau_2^d$};
\draw (9.92,-0.96) node[anchor=north west] {$\tau_2^u$};
\draw (6.27,3.92) node[anchor=north west] {$Y_1$};
\draw (10.19,3.90) node[anchor=north west] {$Y_2$};
\draw (6.13,-1.91) node[anchor=north west] {$\displaystyle M_n=\sum_{\ell=1}^nY_{\ell}$};
\draw [->,dash pattern=on 1pt off 3pt on 6pt off 4pt] (7,4) -- (7,3);
\draw [->,dash pattern=on 1pt off 3pt on 6pt off 4pt] (11,3) -- (11,4);
\draw (4.24,5.61) node[anchor=north west] {$S_n$};
\fill [color=qqqqff] (-2,-2) circle (2.5pt);
\fill [color=qqqqff] (-3,-1) circle (1.5pt);
\fill [color=qqqqff] (-4,0) circle (1.5pt);
\fill [color=qqqqff] (-1,-1) circle (1.5pt);
\fill [color=qqqqff] (0,0) circle (1.5pt);
\fill [color=qqqqff] (-4.26,0.3) circle (0.5pt);
\fill [color=qqqqff] (1,1) circle (1.5pt);
\fill [color=qqqqff] (2,2) circle (1.5pt);
\fill [color=qqqqff] (3,3) circle (1.5pt);
\fill [color=qqqqff] (5,3) circle (1.5pt);
\fill [color=qqqqff] (4,4) circle (2.5pt);
\draw [color=qqqqff] (4,0)-- ++(-1.5pt,0 pt) -- ++(3.0pt,0 pt) ++(-1.5pt,-1.5pt) -- ++(0 pt,3.0pt);
\fill [color=qqqqff] (6,2) circle (2.5pt);
\fill [color=qqqqff] (7,3) circle (2.5pt);
\draw [color=qqqqff] (6,0)-- ++(-1.5pt,0 pt) -- ++(3.0pt,0 pt) ++(-1.5pt,-1.5pt) -- ++(0 pt,3.0pt);
\fill [color=qqqqff] (8,2) circle (1.5pt);
\draw [color=qqqqff] (7,0)-- ++(-1.5pt,0 pt) -- ++(3.0pt,0 pt) ++(-1.5pt,-1.5pt) -- ++(0 pt,3.0pt);
\fill [color=qqqqff] (9,1) circle (2.5pt);
\fill [color=qqqqff] (12,4) circle (2.5pt);
\fill [color=qqqqff] (10,2) circle (1.5pt);
\fill [color=qqqqff] (11,3) circle (1.5pt);
\draw [color=qqqqff] (9,0)-- ++(-1.5pt,0 pt) -- ++(3.0pt,0 pt) ++(-1.5pt,-1.5pt) -- ++(0 pt,3.0pt);
\draw [color=qqqqff] (12,0)-- ++(-1.5pt,0 pt) -- ++(3.0pt,0 pt) ++(-1.5pt,-1.5pt) -- ++(0 pt,3.0pt);
\draw [color=qqqqff] (7,4)-- ++(-1.5pt,0 pt) -- ++(3.0pt,0 pt) ++(-1.5pt,-1.5pt) -- ++(0 pt,3.0pt);
\draw [color=qqqqff] (11,4)-- ++(-1.5pt,0 pt) -- ++(3.0pt,0 pt) ++(-1.5pt,-1.5pt) -- ++(0 pt,3.0pt);
\end{scriptsize}
\end{tikzpicture}
\end{center}
\caption{\label{marche}A one-dimensional PRW}
\end{figure}

Then by a renewal type property (see \cite[Prop. 2.3]{cenac/chauvin/herrmann/vallois/13}), $(\tau_n^d)_{n\geq 1}$ and $(\tau_n^u)_{n\geq 1}$ are independent sequences of {\it i.i.d.\@} random variables. Their distribution tails are straightforwardly given by:
for any $\alpha, \beta\in\{u,d\}$, $\alpha\neq \beta$ and  $n\geq 1$,
\begin{equation}
\label{tail-2d}
\PP(\tau^{\alpha}_{1} \geq n)=\prod_{k=1}^{n-1}q_{\alpha^k\beta}(\alpha).
\end{equation}
Note that Assumption~\ref{ass:a1} amounts to supposing that the \emph{persistence times} $\tau_n^d$ and $\tau_n^u$ are almost surely finite.
The \emph{jump times} (or breaking times) are: $B_0=0$ and, for $n\geq 1$,
\begin{equation}
\label{def:Bn}
B_{2n} := \sum_{k=1}^n \left(\tau_k^d + \tau_k^u \right) \hbox{ and } B_{2n+1} :=  B_{2n} + \tau_{n+1}^d .
\end{equation}
In order to deal with a more tractable random walk built with the possibly unbounded but  \emph{i.i.d.\@} increments $Y_n:=\tau_n^u-\tau_n^d$, we introduce the underlying \emph{skeleton} random walk $(M_n)_{n\geq 1}$ which is the original walk observed at the random times of up-to-down turns:
\begin{equation}
\label{def:skeleton}
M_n:=\sum_{k=1}^{n}Y_k = S_{B_{2n}}.
\end{equation}
Two main quantities play a key role in the asymptotic behaviour, namely the expectations of the lengths of runs:
with Formula \eqref{tail-2d}, let
\begin{equation}
\label{thetaDim1}
\Theta_d:=\E[\tau_1^d] = \sum_{n\geq 1} \prod_{k=1}^{n-1}q_{d^ku}(d) \ \mbox{\ and\ }\ \Theta_u:=\E[\tau_1^u] = \sum_{n \geq 1}\prod_{k=1}^{n-1}q_{u^kd}(u).
\end{equation}
Actually, $\Theta_d$ and $\Theta_u$ already appeared in~\cite[Prop. B1]{cenac/chauvin/herrmann/vallois/13} where it is shown that the driving VLMC of a $1$-dimensional PRW admits a unique invariant probability measure if, and only if $\Theta_d<\infty$ \emph{and} $\Theta_u<\infty$.

Note that  the expectation of $Y_1$ is well defined in $[-\infty,+\infty]$ whenever at least one of the persistence times $\tau_1^u$ or $\tau_1^d$ is integrable.
Thus, as soon as $\Theta_{d}<\infty$ \emph{or} $\Theta_{u}<\infty$, let
\begin{equation}
\label{drift-def1}
{\mathbf d}_{M}:=\E[Y_1] = \underbrace{\Theta_{u}-\Theta_{d}}_{\in[-\infty,+\infty]}
\end{equation}
and
\begin{equation}
\label{drift-def2}
{\mathbf d}_{S}:=\frac{\E[\tau^{u}_{1}]-\E[\tau^{d}_{1}]}{\E[\tau^{u}_{1}]+\E[\tau^{d}_{1}]} =  \frac{\Theta_{u}-\Theta_{d}}{\Theta_{u}+\Theta_{d}}\in [-1,1].
\end{equation}
An elementary computation shows that $\g E\left(M_n\right) = n{\mathbf d}_M$ and $\g E\left(S_n\right) \sim n{\mathbf d}_S$ when $n$ tends to infinity.
Thus, ${\mathbf d}_M$ and ${\mathbf d}_S$ appear as asymptotic drifts when the walks $(M_n)_n$ and $(S_n)_n$ respectively turn out to be transient (see Table~\ref{tableau-rec-trans-2d}).
The behaviour of the walk also depends on quantities $J_{\alpha\mid \beta}$, defined for $\alpha$ and $\beta\in\rond A, \alpha\not=\beta$ by:
\[
J_{\alpha\mid \beta}:=\sum_{n=1}^{\infty}  \frac{n\mathbb \PP(\tau_1^{\alpha}=n)}{\sum_{k=1}^{n} \mathbb \PP(\tau_1^{\beta}\geq k)}.
\]
A complete and usable characterization of the recurrence and the transience of the PRW in terms of the probabilities to persist in the same direction or to switch is given in Proposition \ref{prop:tableau}. 
Its proof relies on a criterion of Erickson (see \cite{Erickson}), applied to the skeleton walk $\left( M_n\right)_n$ which is simpler to deal with because its increments are independent.
\begin{proposition}
\label{prop:tableau}
Under non-nullness assumption and Assumption~\ref{ass:a1}, the random walk $\left( S_n\right)_n$ is recurrent or transient as described in Table \ref{tableau-rec-trans-2d}. 
\end{proposition}

	\begin{table}[h]
		\begin{center}
		\begin{tabular}{|c|c|c|c|c|}
			\hline
			& \multicolumn{2}{c|}{$\Theta_u < \infty$} & \multicolumn{2}{c|}{$\Theta_u = \infty$} \\
			\hline
			\multirow{4}*{$\Theta_d < \infty$} & & drifting $+\infty$ & \multicolumn{2}{c|}{\multirow{4}*{drifting $+\infty$}} \\
			& recurrent & $\mathbf{d}_S > 0$ & \multicolumn{2}{c|}{\multirow{2}*{drifting $+\infty$}} \\
			\cline{3-3}
			& $\mathbf{d}_S=0$ & drifting $-\infty$ & \multicolumn{2}{c|}{}\\ 
			& & $\mathbf{d}_S < 0$ & \multicolumn{2}{c|}{}\\  
			\hline
			\multirow{4}*{$\Theta_d = \infty$} & \multicolumn{2}{c|}{\multirow{4}*{drifting $-\infty$}} &  & drifting $+\infty$ \\
			& \multicolumn{2}{c|}{} & recurrent & $\infty = J_{u\mid d} > J_{d\mid u}$ \\
			\cline{5-5}
			& \multicolumn{2}{c|}{} & $J_{u\mid d}=J_{d\mid u}=\infty$ & drifting $-\infty$ \\
			& \multicolumn{2}{c|}{} & & $\infty = J_{d\mid u} > J_{ u\mid d}$\\
			\hline
		\end{tabular}
\caption{\label{tableau-rec-trans-2d}Recurrence versus Transience (drifting) for $(S_n)_n$ in dimension 1.}
\end{center}
	\end{table}
The most fruitful situation emerges when both running times $\tau_1^u$ and $\tau_1^d$ have infinite means.
In that case, the recurrence properties of $\left( S_n\right)_n$ are related to the behaviour of the skeleton random walk $\left( M_n\right)_n$ defined in \eqref{def:skeleton}, the drift of which, $\mathbf d_{M}$, is not defined.
Thus the behaviour of $\left( S_n\right)_n$ depends on the comparison between the distribution tails of $\tau_1^u$ and $\tau_1^d$ defined in \eqref{tail-2d}, expressed by the quantities~$J_{\alpha\mid \beta}$. Notice that the case when both $J_{u\mid  d}$ and $J_{d\mid  u}$ are finite does not appear in the table since it would imply that $\Theta_u <\infty$ and $\Theta_d <\infty$ (see \cite{Erickson}).

In all three other cases, the drift $\mathbf d_{S}$ is well defined and the PRW is recurrent if and only if $\mathbf d_{S}=0$. In that case, $\displaystyle \lim_{n \to \infty} \frac{S_n}{n}=\mathbf d_{S}=0$. 
Notice that, modifying one transition $q_c$ transforms 
a recurrent PRW into a transient one, since $\mathbf d_{S}$ becomes non-zero.

\subsection{Persistent Random Walks in dimension two}
\label{subsec:PRWdim2}
Take the  alphabet $\mathcal{A}:=\{\ttN,\ttE,\ttW,\ttS\}$. 
Here, $({\mathtt e},{\mathtt n})$ stands for the canonical basis of~$\mathbb Z^2$, ${\mathtt w}=-{\mathtt e}$ and ${\mathtt s}=-{\mathtt n}$. 
Hence, the letters $\ttE$, $\ttN$, $\ttW$ and $\ttS$  stand for moves to the east, north, west and south respectively. Having in mind a random walk with increments in~$\rond A$, any word of the form $\alpha \beta$, $\alpha,\beta\in\rond A, \alpha\not=\beta$ is called a \emph{bend}.
For the sake of simplicity, we condition the walk to start a.s. with a $\ttN\ttE$ bend:
 $\set{X_{-1}=\ttN, X_{0}=\ttE}$.
\begin{figure}[h]
\definecolor{xdxdff}{rgb}{0.49,0.49,1}
\definecolor{ududff}{rgb}{0,0,1}
\begin{center}
\begin{tikzpicture}[scale=0.8,line cap=round,line join=round,>=latex,x=1.0cm,y=1cm]
\clip(-4.3,-4.7) rectangle (14.56,6.34);
\begin{scriptsize}
\draw [->,line width=0.4pt,color=gray] (1.02,-2.26) -- (1,-1);
\draw [line width=0.4pt,color=gray] (1,2)-- (1.02,-1);
\draw [->,line width=2.5pt] (1,1) -- (1,2);
\draw [->,line width=2.5pt] (1,2) -- (2,2);
\draw (1.22,1.78) node[anchor=north west] {$J_0=\ttN\ttE$};
\draw [line width=0.4pt,color=gray] (2,2)-- (5,2);
\draw [->,line width=2.5pt] (5,2) -- (6,2);
\draw [->,line width=2.5pt] (6,2) -- (6,3);
\draw (5.56,1.52) node[anchor=north west] {$J_1=\ttE\ttN$};
\draw [line width=0.4pt,color=gray] (6,4)-- (6,3);
\draw [->,line width=2.5pt] (6,4) -- (6,5);
\draw [->,line width=2.5pt] (6,5) -- (5,5);
\draw [line width=0.4pt,color=gray] (5,5)-- (-2,5);
\draw [->,line width=2.5pt] (-2,5) -- (-3,5);
\draw [->,line width=2.5pt] (-3,5) -- (-3,4);
\draw [line width=0.4pt,color=gray] (-3,4)-- (-3,0);
\draw [->,line width=2.5pt] (-3,0) -- (-3,-1);
\draw [->,line width=2.5pt] (-3,-1) -- (-2,-1);
\draw [line width=0.4pt,color=gray] (-2,-1)-- (3,-1);
\draw [->,line width=2.5pt] (3,-1) -- (4,-1);
\draw [->,line width=2.5pt] (4,-1) -- (4,0);
\draw [line width=0.4pt,color=gray] (4,0)-- (4,3);
\draw [->,line width=2.5pt] (4,3) -- (4,4);
\draw [->,line width=2.5pt] (4.2,4) -- (4.2,3);
\draw [line width=0.4pt,color=gray] (4.2,3)-- (4.2,2);
\draw [->,line width=2.5pt] (4.2,2) -- (4.2,1);
\draw [->,line width=2.5pt] (4.2,1) -- (3.2,1);
\draw [line width=0.4pt,color=gray] (3.2,1)-- (-4.2,1.04);
\draw [->,line width=0.4pt,color=gray] (-3,1.03) -- (-4.2,1.04);
\draw (5.44,5.94) node[anchor=north west] {$J_2=\ttN\ttW$};
\draw (-2.66,4.72) node[anchor=north west] {$J_3=\ttW\ttS$};
\draw (-2.66,-0.32) node[anchor=north west] {$J_4=\ttS\ttE$};
\draw (2.3,-0.36) node[anchor=north west] {$J_5=\ttE\ttN$};
\draw (2.4,3.98) node[anchor=north west] {$J_6=\ttN\ttS$};
\draw (3.94,0.72) node[anchor=north west] {$J_7=\ttS\ttW$};
\draw [->,dash pattern=on 5pt off 5pt] (1,2.5) -- (6,2.5);
\draw (3.34,3) node[anchor=north west] {$B_1$};
\draw [->,dash pattern=on 5pt off 5pt] (6,2.5) -- (1,2.5);
\draw [->,dash pattern=on 5pt off 5pt] (6.5,2) -- (6.5,5);
\draw (6.68,3.92) node[anchor=north west] {$B_2-B_1$};
\draw [->,dash pattern=on 5pt off 5pt] (6.5,5) -- (6.5,2);
\fill [color=ududff] (1,2) circle (2.5pt);
\draw[color=ududff] (0.46,1.96) node {$M_0$};
\fill [color=ududff] (6,2) circle (2.5pt);
\draw[color=ududff] (6.19,1.66) node {$M_1$};
\fill [color=ududff] (6,5) circle (2.5pt);
\draw[color=ududff] (6.32,5.34) node {$M_2$};
\fill [color=ududff] (-3,5) circle (2.5pt);
\draw[color=ududff] (-2.74,5.58) node {$M_3$};
\fill [color=ududff] (-3,-1) circle (2.5pt);
\draw[color=ududff] (-3.12,-1.36) node {$M_4$};
\fill [color=ududff] (4,-1) circle (2.5pt);
\draw[color=ududff] (4.46,-1.2) node {$M_5$};
\fill [color=ududff] (4,4) circle (2.5pt);
\draw[color=ududff] (3.76,4.32) node {$M_6$};
\fill [color=ududff] (4.2,1) circle (2.5pt);
\draw[color=ududff] (4.72,0.92) node {$M_7$};
\end{scriptsize}
\end{tikzpicture}
\end{center}
\vskip -2cm
\caption{\label{marched2}A walk in dimension two.}
\end{figure}
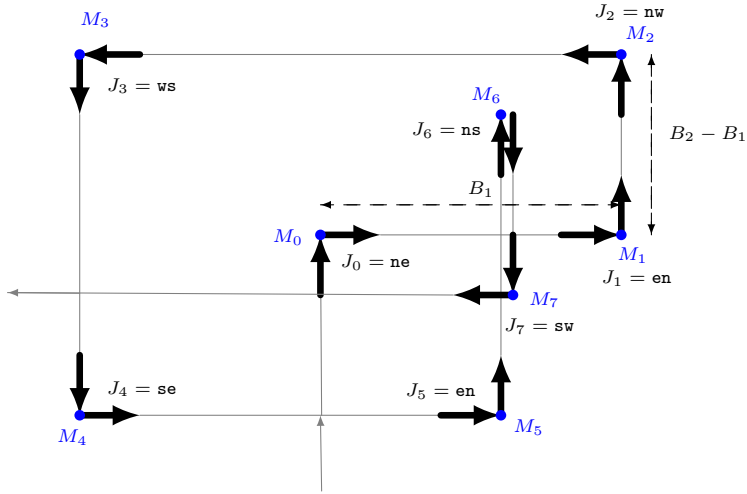

Take a non-null VLMC associated with a quadruple comb on $\rond A$ as drawn in Figure~\ref{fig:bcomb}:
the contexts are $\alpha^n\beta$ for $\alpha,\beta\in\rond A, \alpha\not= \beta$, $n\geq 1$
and the attached probability distributions are denoted by $q_{\alpha^n\beta}$.
The $2$-dimensional PRW $\left( S_n\right) _n$ is defined, using this VLMC, as in Formula \eqref{rw-S}.

Contrary to the $1$-dimensional PRWs, as detailed below, the probability to change direction depends on the time spent in the current direction but also on the previous direction.
As in dimension one, we intend to avoid that $S$ remains frozen in one of the four directions with a positive probability.
Therefore, we make the following assumption, analogous to Assumption~\ref{ass:a1} in dimension $2$.
\begin{assumption}[finiteness of the length of runs]
\label{ass:a2}
For any $\alpha,\beta\in \{\ttN,\ttE,\ttW,\ttS\}$, $\alpha\neq \beta$, 
\begin{equation}
\label{a2}
\lim _{n\to +\infty}\left(\prod_{k=1}^{n}q_{\alpha^k\beta}(\alpha)\right)=0.
\end{equation}
\end{assumption}
Let $(B_n)_{n\geq 0}$ be the \emph{breaking times} defined inductively by 
\begin{equation}\label{def:jump}
B_0=0\quad\mbox{and}\quad
B_{n+1}=\inf\left\{k>B_{n} : X_{k}\neq X_{k-1}\right\}.
\end{equation}
As in dimension $1$, Assumption~\ref{ass:a2} implies that the breaking times $B_n$ are almost surely finite.

Define the so called \emph{internal chain} $\left( J_n\right) _{n\geq 0}$ by $J_0=\ttN\ttE$ and, for all $n\geq 1$,
\begin{equation}\label{drivingchain}
J_{n}:=X_{B_{n-1}}X_{B_{n}}.
\end{equation}
Let us illustrate these random variables by a small example, in which:
$B_1=4$, $B_2=7$, $J_0=X_{-1}X_0$, $J_1=X_{B_0}X_{B_1}=X_0X_4$, $J_2=X_{B_1}X_{B_2}=X_4X_7$.

\begin{tikzpicture}
\newcommand{\hs}{0.7}
\newcommand{\vs}{-0.5}
\draw (-1.2*\hs,0) node{$-1$};
\foreach \j in {0,1,2,3,4,5,6,7} \draw (\j*\hs,0) node{$\j$};
\draw (-1.2*\hs,\vs) node{$\ttN$};
\foreach \j in {0,1,2,3} \draw (\j*\hs,\vs) node{$\ttE$};
\foreach \j in {4,5,6} \draw (\j*\hs,\vs) node{$\ttN$};
\foreach \j in {7} \draw (\j*\hs,\vs) node{$\ttW$};
\foreach \j in {0,1,2,3} \draw (\j*\hs,2*\vs) node{$\ttN\ttE$};
\foreach \j in {4,5,6} \draw (\j*\hs,2*\vs) node{$\ttE\ttN$};
\foreach \j in {7} \draw (\j*\hs,2*\vs) node{$\ttN\ttW$};
\draw (-0.5*\hs,0)--++(0,-2);
\draw (3.5*\hs,0)--++(0,-2);
\draw (6.5*\hs,0)--++(0,-2);
\draw (0.2,3.5*\vs) node{\small$B_0=0$};
\draw (0.3,4.5*\vs) node{\small$J_0=\ttN\ttE$};
\draw (4*\hs+0.2,3.5*\vs) node{\small$B_1=4$};
\draw (4*\hs+0.3,4.5*\vs) node{\small$J_1=\ttE\ttN$};
\draw (7*\hs+0.2,3.5*\vs) node{\small$B_2=7$};
\draw (7*\hs+0.3,4.5*\vs) node{\small$J_2=\ttN\ttW$};
\draw (-2.5*\hs,0) node{$n$:};
\draw (-2.5*\hs,\vs) node{$X_n$:};
\draw (-2.5*\hs,2*\vs) node{$Z_n$:};
\end{tikzpicture}

The process $\left( J_n\right) _{n\geq 0}$ is an irreducible Markov chain on the set of bends $\mathcal{S}:=\set{\alpha\beta|\alpha \in \CA, \beta \in \CA, \alpha\neq \beta}$.
Its Markov kernel is defined by: for every $\beta,\alpha,\gamma\in\mathcal A$ with $\beta\neq\alpha$ and $\alpha\neq\gamma$,
\begin{equation}
\label{markovsymb1}
P(\beta\alpha;\alpha\gamma):=\sum_{n=1}^{\infty}\left(\prod_{k=1}^{n-1}q_{\alpha^k \beta}(\alpha)\right)q_{\alpha^n \beta}(\gamma),
\end{equation}
the numbers $P(\alpha\beta,\gamma\delta)$ being $0$ for every couple of bends not of the previous form.
Remark that the non-nullness assumption (see Definition~\ref{nn}) implies the irreducibility of $\left( J_n\right)_n$ and its aperiodicity.
The state space $\CS$ is finite so that $\left( J_n\right)_n$ is positive recurrent: it admits a unique invariant probability measure $\pi_{\scriptscriptstyle J}$. 

Denote $T_0=0$ and $T_{n+1}:=B_{n+1}- B_n$ for every $n\geq 0$.
These waiting times (also called \emph{persistence times}) are not independent, contrary to the one-dimensional case.
The \emph{skeleton random walk} $(M_{n})_{n\geq 0}$ on $\mathbb Z^2$ -- which is the  PRW observed at the breaking times --  is then defined as 
\begin{equation}
M_n:=S_{B_n}=\sum_{i=1}^{n}\left(\sum_{k=B_{i-1}+1}^{B_{i}} X_k\right) = \sum_{i=1}^{n}\left( B_{i} - B_{i-1}\right) X_{B_i} .
\end{equation}
Notice that $(M_n)_n$ is generally not a classical RW with \emph{i.i.d.} increments. Nevertheless, taking into account the additional information given by the internal Markov chain $(J_n)_n$, then $\left( J_n,M_n\right)_n$ is a Markov Additive Process (see \cite{cinlar/72})
as it will appear in Section~\ref{sec:meetic}.

Here, $(J_n)_n$ is positive recurrent but this
does not imply the recurrence of $(S_n)_n$ or $(M_n)_n$.
Moreover, $(S_n)_n$ and $(M_n)_n$ may have different behaviours.
Explicit necessary and sufficient conditions for the recurrence of $(M_n)_n$ in terms of characteristic functions and convergence of suitable series are given in \cite[Theorem 2.1]{cenac:hal-01658494}.
The following proposition states a dichotomy between some recurrence \emph{versus} transience phenomenon. 
\begin{theorem}
Under non-nullness assumption, the following dichotomy holds.

(i) The series $\sum _n\mathbb P\left( M_n=0\right)$ diverges if, and only if the process $\left( M_n\right) _n$ is recurrent in the following sense:
\[
\exists r>0,~~\PP\left(\liminf_{n\to\infty} \|M_n\|<r\right)=1.
\]

(ii) The series $\sum _n\mathbb P\left( M_n=0\right)$ converges if, and only if the process $\left( M_n\right) _n$ is transient in the following sense:
\[
\PP\left(\lim_{n\to\infty} \|M_{n}\|=\infty\right)=1.
\]
\end{theorem}

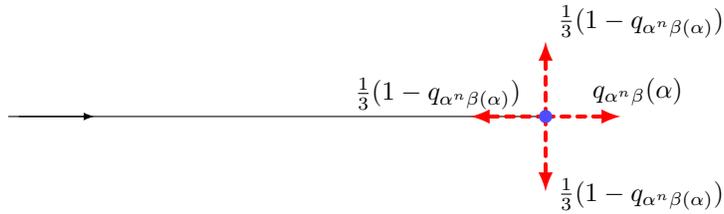
\begin{figure}[h]
\vskip -2cm
\definecolor{ffqqqq}{rgb}{1,0,0}
\definecolor{ududff}{rgb}{0.3,0.3,1}
\begin{tikzpicture}[line cap=round,line join=round,>=latex,x=1.0cm,y=1.0cm]
\clip(-4.3,-3.28) rectangle (8.68,6.3);
\draw (-4.14,2)-- (3,2);
\draw [->] (-4,2) -- (-3,2);
\draw (3.5,2.64) node[anchor=north west] {$q_{\alpha^n\beta}(\alpha)$};
\draw (3.04,3.6) node[anchor=north west] {$\frac{1}{3}(1-q_{\alpha^n\beta(\alpha)})$};
\draw (3.04,1.3) node[anchor=north west] {$\frac{1}{3}(1-q_{\alpha^n\beta(\alpha)})$};
\draw (0.34,2.64) node[anchor=north west] {$\frac{1}{3}(1-q_{\alpha^n\beta(\alpha)})$};
\draw [->,line width=1.5pt,dash pattern=on 3pt off 3pt,color=ffqqqq] (3,2) -- (3,3);
\draw [->,line width=1.5pt,dash pattern=on 3pt off 3pt,color=ffqqqq] (3,2) -- (4,2);
\draw [->,line width=1.5pt,dash pattern=on 3pt off 3pt,color=ffqqqq] (3,2) -- (3,1);
\draw [->,line width=1.5pt,dash pattern=on 3pt off 3pt,color=ffqqqq] (3,2) -- (2,2);
\begin{scriptsize}
\fill [color=ududff] (3,2) circle (2.5pt);
\end{scriptsize}
\end{tikzpicture}
\vskip -4cm
\caption{\label{transitions2}The original Directionally Reinforced Random Walk (DRRW).}
\end{figure}

Does the recurrence (resp. the transience) of $(M_n)_n$ and $(S_n)_n$ occur at the same time?
The answer to this twenty-year-old question is no.
\begin{theorem}
[Definitive invalidation of the conjecture in \cite{Mauldin1996}]
\label{conjecture}
	There exist recurrent PRWs  $(S_n)_n$ 
	having an associated transient MRW skeleton $(M_n)_n$.
\end{theorem}
Supposing that the persistence time distributions are horizontally and vertically symmetric is a natural necessary condition for the random walk $(S_n)_n$ to be recurrent.
One example is given by the Directionally Reinforced Random Walk (DRRW), originally introduced in \cite{Mauldin1996}, see Figure~\ref{transitions2}.
Some particular values of the transition probabilities $q_{\alpha ^n\beta}$ provide counterexamples.
It is shown in~\cite{cenac:hal-01658494} that the corresponding distributions of the persistence times must be
non-integrable.
In Section~\ref{sec:meetic}, this non integrability will be related to non existence of any invariant probability measure for the driving VLMC.




\section{VLMC: existence of stationary probability measures}
\label{sec:VLMC}

Take a VLMC denoted by $U=\left( U_n\right) _{n\geq 0}$, defined by a pair $\left(\CT,q\right)$ where $\CT$ is a context tree on an alphabet~$\CA$ and $q=\left( q_c\right) _{c\in\CC}$ a family of probability measures on $\CA$, indexed by the contexts of $\CT$.
A probability measure $\pi$ on $\CR$ is \emph{stationary} or \emph{invariant} (with regard to $U$) whenever $\pi$ is the distribution of every $U_n$ as soon as it is the distribution of $U_0$.
The question of interest consists here  in finding conditions on $\left(\CT,q\right)$ for the process to admit at least one -- or a unique one -- stationary probability measure.
The heuristic presentation aims to show how combinatoric objects -- namely the $\alpha$-lis of contexts -- and numbers -- the cascades -- naturally emerge.

Assume that $\pi$ is a stationary probability measure on $\CR$.

$\bullet$ First step: finite words.
Since $\CR$ is endowed with the cylinder $\sigma$-algebra, $\pi$ is determined by its values $\pi\left( w\CR\right)$ on the cylinders $w\CR$, where $w$ runs over all finite words on $\CA$.

$\bullet$ Second step: longest internal suffixes of words.
Assume that $e$ is a finite non-internal word and take $a\in\CA$.
Then, its $\pref$ is well defined and, because of Formula~\eqref{probaTransition}, since $\pi$ is stationary,
\begin{equation}
\label{preCascade}
\pi\left(\alpha e \CR\right)=q_{\pref(e)}(\alpha )\times\pi\left(e\CR\right).
\end{equation}
Iterating this formula as far as possible leads to the following definitions.
Consider any non-empty finite word $w$.
It is uniquely decomposed as $w=p\alpha s=\beta _1\beta_2\beta _3\cdots\beta_\ell\alpha s$, where $\alpha$ and the $\beta$'s are letters and $s$ is the \emph{longest internal suffix} of $w$.
The integer $\ell$ is non-negative and $p=\beta _1\beta_2\cdots\beta_\ell$ is a prefix of $w$ that may be empty -- in which case $\ell =0$.

\begin{definition}
[lis and $\alpha$-lis]
With these notations, the Longest Internal Suffix $s$ is shortened as the $\emph{lis}$ of $w$.
The word $\alpha s$ is called the \emph{$\alpha$-lis} of $w$.
\end{definition}

\begin{definition}
[cascade]
With the notation above, the \emph{cascade} of $w$ is the product
\begin{equation}
\label{formCascade}
\casc (w)
=q_{\pref(\beta_2\cdots\beta_\ell\alpha s)}(\beta _1 )
q_{\pref(\beta_3\cdots\beta_\ell\alpha s)}(\beta _2 )
\cdots
q_{\pref(\alpha s)}(\beta_\ell).
\end{equation}
\end{definition}
Note that this definition makes sense because all the $\beta_k\cdots\beta_\ell\alpha s$ are non-internal words, $k\geq 2$.
Moreover, if $w=\alpha s$ where $s$ is internal, then $\ell =0$ and $\casc (w)=1$.
With these definitions, iterating Formula~\eqref{preCascade} leads to the following equality, named \emph{Cascade Formula}:
for every non-empty finite word $w$ having $\alpha s$ as an $\alpha$-lis,
\begin{equation}
\pi\left(w \CR\right) = \casc (w)\times\pi\left(\alpha s\CR\right).
\end{equation}
This shows that $\pi$ is determined by its values on words of the form $\alpha s$ where $s$ is internal and $\alpha\in\CA$.

$\bullet$ Third step: finite contexts.
Assume that $s$ is an internal word and that $\alpha\in\CA$.
It is shown in~\cite{cenac/chauvin/paccaut/pouyanne/18} that a stationary probability measure never charges infinite words so that, by disjoint union,
\begin{equation}
\label{chaisPasCommentLApeller}
\pi\left(\alpha s \CR\right)=\sum _{\substack{
		{c:}{\rm~finite~context}\\
		c=s\cdots}}
\pi\left(\alpha c\CR\right)=\sum _{\substack{
		{c:}{\rm~finite~context}\\
		c=s\cdots}}
q_c(\alpha )\pi\left(c\CR\right).
\end{equation}
Note that the set of indices may be infinite but the family is summable because $\pi$ is a finite measure.
This shows that $\pi$ is entirely determined by its values $\pi\left(c\CR\right)$ on the 
finite contexts.

$\bullet$ Fourth step: $\alpha$-lis of finite contexts.
Cascade Formula~\eqref{formCascade} applied to any finite context $c$
(contexts are non-empty words) writes $\pi\left(c\CR\right)=\casc (c)\pi\left(\alpha _cs_c\CR\right)$, where $\alpha _cs_c$ is the $\alpha$-lis of~$c$.
Denote by $\CS =\CS\left(\CT\right)$ the set of finite context $\alpha$-lis:
\[
\CS=\set{\alpha _cs_c:~c{\rm ~finite~context}}.
\]
If $s$ is an internal word and if $\alpha\in\CA$, then Formula~\eqref{chaisPasCommentLApeller} leads to
\begin{equation}
\label{rhoooYenAMarreDeChercherDesNomsDeFormules}
\pi\left(\alpha s \CR\right)
=\sum _{\substack{
		{c:}{\rm~finite~context}\\
		c=s\cdots}}
\casc\left(\alpha c\right)\pi\left(\alpha _cs_c\CR\right),
\end{equation}
showing that $\pi$ is determined by its values $\pi\left(\alpha _cs_c\CR\right)$ on $\CS$.

$\bullet$ Last step: a (generally infinite) linear system.
When $w$ and $v$ are finite words and when $\alpha s\in\CS$, the notation
\[
w=v\cdots =\cdots [\alpha s]
\]
stands for: $w$ has $v$ as a prefix and $\alpha s$ as an $\alpha$-lis.
Writing Formula~\eqref{rhoooYenAMarreDeChercherDesNomsDeFormules} for every $\alpha s\in\CS$ and grouping in each of them the terms that arise from contexts having the same $\alpha$-lis leads to the following square system
(at most countably many unknowns $\pi\left(\alpha s\CR\right)$ and as many equations):
\begin{equation}
\label{systQ}
\forall \alpha s\in\CS,~
\pi\left(\alpha s\CR\right)
=\sum _{\beta t\in\CS}\pi\left(\beta t\CR\right)
\left(\sum _{\substack{{c:}{\rm~finite~context}\\c=s\cdots =\cdots [\beta t]}}\casc\left(\alpha c\right)\right).
\end{equation}

\begin{definition}[Matrix Q]
When $\CT$ is a context tree having $\CS$ as a context $\alpha$-lis set, $Q=Q\left(\CT\right)$ is the $\CS$-indexed square matrix defined by:
\begin{equation}
\label{defQ}
\forall\alpha s,\beta t\in\CS,~
Q_{\beta t,\alpha s}
=\sum _{\substack{{c:}{\rm~finite~context}\\c=s\cdots =\cdots [\beta t]}}\casc\left(\alpha c\right)
\in [0,+\infty ].
\end{equation}
\end{definition}

Thus, System~\eqref{systQ} tells us that, when $\pi$ is a stationary measure, the row-vector $\left(\pi\left(\alpha s \CR\right)\right) _{\alpha s\in\CS}$ appears as a left-fixed vector of the matrix $Q$.

\begin{definition}[Cascade series]
\label{def:cascseries}
For every $\alpha s\in\CS$, denote
\[
\kappa _{\alpha s} =\sum _{\substack{{c:}{\rm~finite~context}\\c=\cdots [\alpha s]}}\casc (c)\in [0,+\infty ].
\]
When this series is summable, one says that \emph{the cascade series of $\alpha s$ converges}.
Whenever the cascades series of all $\alpha s\in\CS$ converge, one says that \emph{the cascade series (of the VLMC) converge}.
\end{definition}
Note that the convergence of (all) the cascade series is sufficient to guarantee the finiteness of $Q$'s entries.
Actually, for a general VLMC, as it is made precise in~\cite{cenac/chauvin/paccaut/pouyanne/18}, the convergence of the cascade series appears as a pivot condition when dealing with existence and unicity of a stationary probability measure.
In this paper, we just state a necessary and sufficient condition for a special kind of VLMC:
the \emph{stable} ones that have a finite $\CS$.
The following proposition is proven in~\cite{cenac/chauvin/paccaut/pouyanne/18}.

\begin{proposition}\label{prop:defstable}
Let $\rond T$ be a context tree. The following conditions are equivalent.
\begin{enumerate}
	\item[(i)] $\forall \alpha \in \rond A$, $\forall w \in \rond W$, $\alpha w \in \rond T \Longrightarrow w \in \rond T$.
		\item[(ii)] If $c$ is a finite context and $\alpha \in \rond A$, then $\alpha c$ is non-internal.
		\item[(iii)] $\rond T\subseteq\CA\CT=\{\alpha w,~\alpha\in\CA,~w\in\CT\}$.
		\item[(iv)] For any VLMC $(U_n)_n$ associated with $\rond T$, the process $\left(\pref(U_n)\right)_{n\in\g N}$ is a Markov chain that has the set of contexts as a state space.
\end{enumerate} 
\end{proposition}

The context tree is called \emph{stable} whenever one of these conditions is fulfilled.

It turns out that the stability of $\rond T$ together with the non-nullness of the VLMC imply both stochasticity and irreducibility of the matrix $Q$.
Consequently, in the simple case where $Q$ is a finite-dimensional matrix, there exists (thanks to stochasticity) a unique (thanks to irreducibility) left-fixed vector for $Q$.
As a consequence of a much more general result proven in~\cite{cenac/chauvin/paccaut/pouyanne/18}, this implies existence and unicity of a stationary probability measure for the VLMC, as stated below.

\begin{theorem}
\label{th:stablefini}
	Let $({\CT},q)$ be a non-null stable probabilized context tree. If $\#{\CS}<\infty$, then the following are equivalent.
\begin{enumerate}
	\item The VLMC associated to $({\CT},q)$ has a unique stationary probability measure.
	\item The cascade series converge (see Definition \ref{def:cascseries}). 
\end{enumerate}
\end{theorem}

Notice that in the non stable case, the matrix $Q$ is generally not stochastic nor is it even substochastic.
Notice also that, even in the stable case, when $\#\rond S=\infty$, the matrix $Q$ may be stochastic, irreducible and positive recurrent while the VLMC does not admit any stationary probability measure.
One can find such an example in~\cite{cenac/chauvin/paccaut/pouyanne/18}, built with a left comb of left comb -- see Example~\ref{ex:4VLMC}.

\section{Where VLMC and PRW meet}
\label{sec:meetic}

On one hand, a VLMC is defined by its context tree and its transition probability distributions $q_c$
-- in particular the double and the quadruple combs which are stable trees with finitely many context $\alpha$-lis.

Necessary and sufficient conditions of existence and uniqueness of stationary probability measures are given in terms of cascade series. On the other hand, for PRW (defined from VLMC), recurrence properties are written in terms of persistence times. Our aim is to build a bridge between these two families of objects and properties. The meeting point turns out to be the semi-Markov processes of $\alpha$-lis and bends.

\subsection{Semi-Markov chains and Markov Additive Processes}
\label{subsec:defSM}
Semi-Markov chains are defined following \cite{barbu/limnios/08} thanks to so-called Markov renewal chains.
\begin{definition}[Markov renewal chain]
\label{def:MRC}
A Markov chain $(J_n,T_n)_{n\geq0}$ with state space $\rond E\times \g N$ is called a (homogeneous) \emph{Markov renewal chain} (shortly MRC) whenever the transition probabilities satisfy: $\forall n\in\g N$, $\forall a,b\in\rond E$, $\forall j,k\in\g N$,
\[
\begin{array}{rl}
\PP\left(  J_{n+1}=b, T_{n+1} = k\big| J_n = a, T_n =j   \right) &= \PP\left(  J_{n+1}=b, T_{n+1} = k\big| J_n = a  \right) \\
[5pt]
&=: p_{a,b}(k)
\end{array}
\]
and $\forall a,b\in\rond E$, $p_{a,b}(0) = 0$.
For such a chain, the family $p=\left(p_{a,b}(k)\right)_{a,b\in\rond A, k\geq 1}$ is called its \emph{semi-Markov kernel}.
\end{definition}

\begin{definition}[Semi-Markov chain]
\label{def:semiMarkov}
Let $(J_n,T_n)_{n\geq0}$ be a Markov renewal chain with state space $\rond E\times \g N$.
Assume that $T_0=0$.
For any $n\in \g N$,  let $B_n$ be defined by
\[
B_n = \sum_{i=0}^n T_i.
\]
The \emph{semi-Markov chain} associated with $(J_n,T_n)_{n\geq0}$ is the $\rond E$-valued process $(Z_j)_{j\geq0}$ defined by 
\[
\forall j \hbox{ such that } B_n\leq j < B_{n+1}, \hskip 5mm Z_j = J_n.
\]
\end{definition}
Note that the sequence $(B_n)_{n\geq 0}$ is almost surely increasing because of the assumption $p_{a,b}(0) = 0$ (instantaneous transitions are not allowed) that guarantees that $T_n\geq 1$ almost surely, for any $n\geq 1$.

The $B_n$ are \emph{jump times}, the $T_n$ are \emph{sojourn times} in a given state and $Z_j$ stagnates at a same state between two successive jump times. The process $\left( J_n\right) _n$ is called the \emph{internal (underlying) chain} of the semi-Markov chain $\left( Z_n\right) _n$.

The previous definitions make transitions to the same state between time $n$ and time $n+1$ possible. Nevertheless, one can boil down to the case where $p_{a,a}(k)=0$ for all $a\in\rond E, k\in\g N$ (see the details in \cite{cenac/chauvin/paccaut/pouyanne/18}). 

A close notion, \emph{Markov Additive Processes}, can be found in \cite{cinlar/72}.
\begin{definition}[Markov Additive Process]
\label{def:MAP}
A Markov chain $(J_n,B_n)_{n\geq0}$ with state space $\rond E\times \g N$ is called a \emph{Markov Additive Process} (shortly MAP) whenever  $\left( J_n,B_{n}-B_{n-1}\right) _n$ is a 
Markov renewal chain.
\end{definition}
\subsection{Persistent Random Walks induce semi-Markov chains}
\label{subsec:PRWandSM}
Let us start with $1$-dimensional PRW, as defined in Section \ref{subsec:PRWdim1}.
In this case, at each time $j, j\geq 0$, the increment $X_j$ of the walk $S$ takes $d$ or $u$ as a value (see Figure \ref{marche}).
Let us see that $(X_j)_{j\geq 0}$ is a semi-Markov chain, starting from $X_0 = d$. 
Remember that $B_n$ denotes the $n$-th jump times -- see Equation \eqref{def:Bn}. 
Define then $(J_n)_n$ by 
\begin{equation}
\label{def:JnDim1}
J_n:=X_{B_n}.
\end{equation}
Moreover, let $T_n$ be the $n$-th waiting time, namely $T_0=0$ and, for $n\geq 1$,
\[
T_n=B_n-B_{n-1}.
\]
These waiting times are related to the persistence times $\tau$ by the following formulae: for all $k\geq 1$,
\begin{equation}
\label{def:TnDim1}
T_{2k}:=\tau_k^u \hbox{  \ \  and  \ \   } T_{2k-1}:=\tau_k^d.
\end{equation}
With these notations, $(J_n,T_n)_{n\geq 0}$ is a Markov renewal chain and its semi-Markov kernel writes: $\forall \alpha,\beta\in\{ u,d\},\alpha\not= \beta,\forall k\geq 1$,
\begin{equation}
\label{SMkernelDim1}
p_{\alpha,\beta}(k) =  \left(\prod_{j=1}^{k-1}q_{\alpha^j \beta}(\alpha)\right)q_{\alpha^k \beta}(\beta),
\end{equation}
as can be straightforwardly checked.
Moreover, Assumption~\ref{ass:a1} guarantees that the $T_n$ are a.s. finite.
Besides, Formulae~\eqref{thetaDim1} write
\[
\g E\left( T_{2k}\right)=\Theta _u
{\rm ~and~~}
\g E\left( T_{2k+1}\right)=\Theta _d.
\]
The situation in dimension $1$ is summarized by the following proposition.
\begin{proposition}
\label{prop:semiMarkovPRW1}
For a PRW in dimension 1, defined by a VLMC associated with a double comb, the sequence $(X_j)_j$ of the increments is an  $\rond A$-valued semi-Markov chain with Markov renewal chain $(J_n,T_n)_n$ as defined in \eqref{def:JnDim1} and \eqref{def:TnDim1} and its semi-Markov kernel is given by equation \eqref{SMkernelDim1}.
\end{proposition}

Let us deal now with the $2$-dimensional PRW, defined in Section~\ref{subsec:PRWdim2}. At each time~$j, j\geq 0$, the increment $X_j$ of the walk $S$ takes $\ttN, \ttE, \ttW$ or $\ttS$ as a value.
But, as already noticed, changing direction depends on the time spent in the current direction but also, contrary to the $1$-dimensional PRWs, on the previous direction.
In otherwords, the bends play the main role.
This gives rise to the process $\left( Z_j\right) _j$, valued  in the set of bends $\{\alpha\beta : \alpha, \beta\in\rond A, \alpha\not=\beta \}$, defined in the following manner:
$Z_0 = X_{-1}X_0=\ttN\ttE$ and, for $j\geq 1$, $Z_j = \alpha\beta$ if and only if $X_j = \beta$
and the first letter distinct from $\beta$ in the sequence $X_{j-1},X_{j-2},X_{j-3},\cdots$ is $\alpha$.
Let us see that $(Z_j)_{j\geq 0}$ is a semi-Markov chain.
Use here notations $(J_n)_n$, $(B_n)_n$ and $(T_n)_n$ of section~\ref{subsec:PRWdim2}.

Notice that, contrary to the one-dimensional case, the waiting times $T_n$ are not independent.
Nevertheless, $(J_n,T_n)_{n\geq 0}$ is a Markov renewal chain with semi-Markov kernel
\begin{equation}
\label{SMkernelDim2}
p_{\beta\alpha,\alpha\gamma}(k) := \left(\prod_{j=1}^{k-1}q_{\alpha^j \beta}(\alpha)\right)q_{\alpha^k \beta}(\gamma),
\end{equation}
as can be straightforwardly checked. Summarizing, the following proposition holds.

\begin{proposition}
\label{prop:semiMarkovPRW2}
For a PRW in dimension 2, defined by a VLMC associated with a quadruple comb, the sequence $(Z_j)_j$ of the bends is a semi-Markov chain with Markov renewal chain $(J_n,T_n)_n$ as defined in Section~\ref{subsec:PRWdim2}.
Its semi-Markov kernel is given by equation \eqref{SMkernelDim2}. In addition, $(J_n,B_n)_n$ is a Markov Additive Process.
\end{proposition}

\subsection{Semi-Markov chain of the $\alpha$-lis in a stable VLMC}
\label{subsec:compareSM}

In this section, let us consider a more general case than a double comb or a quadruple comb, namely a stable VLMC. In this case, there is always a semi-Markov chain induced by the process $(U_n)_n$, as described in the following.

Let $(U_n)_{n\geq 0}$ be a stable non-null VLMC such that the series of cascades converge (see Definition \ref{def:cascseries}).
Recall that $\rond S$ denotes the set of context $\alpha$-lis of the VLMC.
Let $(C_n)_{n\geq 0}$ be the sequence of contexts and for $n\geq 0$, let ${Z}_n$ be the  $\alpha$-lis of $C_n$:
\[
C_n = \pref (U_n)  \hskip 5mm \hbox{ and }  \hskip 5mm {Z}_n = \alpha_{C_n}s_{C_n}.
\]
\begin{proposition}
\label{pro:compareSM}
Let $({B}_n)_{n\geq 0}$ be the increasing sequence of times defined by ${B}_0 = 0$ and for any $n\geq 1$,
\[
{B}_n = \inf \set{ k> {B}_{n-1}, |C_k| \leq |C_{k-1}|}=  \inf \set{ k> {B}_{n-1}, C_k \in\rond S}
\]
and let ${T}_n = {B}_n - {B}_{n-1}$ for $n\geq 1$ and ${T}_0 = 0$. For any $n\geq 0$, let ${J}_n = {Z}_{{B}_n}$. Then
\begin{itemize}
\item[(i)]
 ${B}_n$ and $T_n$ are almost surely finite and for $\alpha s\in\rond S$, $\g E\left( T_n \big| J_n = \alpha s\right) = \kappa_{\alpha s}$.
\item[(ii)]
$({Z}_n)_{n\geq 0}$ is an $\rond S$-valued semi-Markov chain associated with the Markov renewal  chain $({J}_n, {T}_n)_{n\geq 0}$.
\item[(iii)]
The associated semi-Markov kernel writes: $\forall \alpha s, \beta t \in\rond S$, $\forall k\geq 1$,
\[
p_{\alpha s, \beta t}(k) = \sum_{\substack{c\in\rond C,~c=t\cdots\\c= \cdots [\alpha s]\\|c| = |\alpha s| + k-1}} \casc\left(\beta c\right) .
\]
\end{itemize}
\end{proposition}

The proof is detailed in~\cite{cenac/chauvin/paccaut/pouyanne/18}.
It relies on the way the VLMC grows between two jump times: at the beginning, letters are added to the current context $C_n$, the $\alpha$-lis does not change and the length of the current context increases one by one. At a certain time (a.s. finite), adding a letter to the current context does not provide a context any more but an external node. At this moment, it happens (it is not trivial and only holds for a stable context tree) that

(i) the $\alpha$-lis of the current context is renewed;

(ii) the length of current context does not grow;

(iii) the current context begins by a lis.

These mechanisms explain the expressions of $B_n$ and the formula giving the semi-Markov kernel.

\begin{remark}
\label{rem:semiMarkov}
In the very particular case of the double or quadruple comb, the semi-Markov chain $\left( Z_n\right) _n$ contains as much information as the chain~$\left( U_n\right) _n$.
But in general, the semi-Markov chain $\left( Z_n\right) _n$ contains less information than the chain $\left( U_n\right) _n$. 
To illustrate this, here is an example with a finite context tree.

\vskip 5pt
\begin{minipage}{0.3\textwidth}
\centering
\begin{tikzpicture}[scale=0.4]
		\tikzset{every leaf node/.style={draw,circle,fill},every internal node/.style={draw,circle,scale=0.01}}
\Tree [.{} [.{} [.{} {} [.{} {} {} ] ] [.{} {} [.{} {} {} ] ] ] [.{} {} [.{} {} {} ] ] 
	]
\end{tikzpicture}
\end{minipage}
\begin{minipage}{0.6\textwidth}
\centering
\begin{tabular}{r|l}
$\alpha$-lis $\alpha s$&contexts having $\alpha s$ as an $\alpha$-lis  \\
\hline
10&10,010,110,0010,0110\\
000&000\\
111&111,0111\\
0011&0011
\end{tabular}
\end{minipage}
\vskip 5pt
In this example, 0010 and 0110 are two contexts of the same length, with the same $\alpha$-lis 10 and beginning by the same lis 0. Hence if we know that ${J}_n=10$, ${B}_{n+1}-{B}_{n}=3$ and ${J}_{n+1}=10$, then ${Z}_j$ is uniquely determined between the two successive jump times, whereas
	there are two possibilities to reconstruct the VLMC $(U_n)_n$. With the notations of Proposition \ref{pro:compareSM}, there are two cascade terms in $p_{10,10}(3) $:

\begin{minipage}{0.68\textwidth}
\begin{align*}
p_{10,10}(3) &= \PP\left(C_{{B}_n+1} = 010, C_{{B}_n+2} = 0010,C_{{B}_n+3} = 10010 | C_{{B}_n} =10\right)\\
& \ \ + \PP\left(C_{{B}_n+1} = 110, C_{{B}_n+2} = 0110,C_{{B}_n+3} = 10110 | C_{{B}_n} =10\right)\\
&= q_{10}(0) q_{010}(0) q_{0010}(1) + q_{10}(1) q_{110}(0) q_{0110}(1)\\
&= \casc (10010) + \casc (10110).
\end{align*}
\end{minipage}
\end{remark}

\subsection{The meeting point}
\label{subsec:meetic}

Summing up, the announced close encounter can be done with the following (commutative) diagram, together with the following explanations.

\begin{equation}
\begin{tikzpicture}
\label{diagrammeComm}
\newcommand{\largeur}{4.6}
\newcommand{\hauteur}{1.5}
\newcommand{\lgFlecheHoriz}{1}
\newcommand{\lgFlecheVert}{0.5}
\draw (0,\hauteur) node{MRC $\left( {J_n^V},T_n\right)_n$};
\draw (0,2*\hauteur) node{VLMC $\left( U_n\right)_n$};
\draw (\largeur,\hauteur) node{MAP $\left( {J_n^W},M_n\right)_n$};
\draw (\largeur,2*\hauteur) node{PRW $\left( S_n\right)_n$};
\draw (\largeur,0) node{Semi-Markov $\left( Z_n^W\right)_n$};
\draw (0,0) node{Semi-Markov $\left( Z_n^V\right)_n$};
\draw (0.5*\largeur,2*\hauteur)--++(-0.5*\lgFlecheHoriz,0);
\draw [->,>=latex] (0.5*\largeur,2*\hauteur)--++(0.5*\lgFlecheHoriz,0);
\draw (0.5*\largeur,2*\hauteur+0.3) node{$D$};
\draw (0.5*\largeur,\hauteur)--++(0.5*\lgFlecheHoriz,0);
\draw [->,>=latex] (0.5*\largeur,\hauteur)--++(-0.5*\lgFlecheHoriz,0);
\draw (0.5*\largeur,\hauteur+0.3) node{$N$};
\draw [->,>=latex] (0.5*\largeur,0)--++(-0.5*\lgFlecheHoriz,0);
\draw [->,>=latex] (0.5*\largeur,0)--++(0.5*\lgFlecheHoriz,0);
\draw (0.5*\largeur,0+0.3) node{$R$};
\draw (0,1.5*\hauteur)--++(0,0.5*\lgFlecheVert);
\draw [->,>=latex] (0,1.5*\hauteur)--++(0,-0.5*\lgFlecheVert);
\draw (0.3,1.5*\hauteur) node{$L$};
\draw (0,0.5*\hauteur)--++(0,0.5*\lgFlecheVert);
\draw [->,>=latex] (0,0.5*\hauteur)--++(0,-0.5*\lgFlecheVert);
\draw (0.4,0.5*\hauteur) node{$S_V$};
\draw (\largeur,1.5*\hauteur)--++(0,0.5*\lgFlecheVert);
\draw [->,>=latex] (\largeur,1.5*\hauteur)--++(0,-0.5*\lgFlecheVert);
\draw (\largeur+0.3,1.5*\hauteur) node{$B$};
\draw (\largeur,0.5*\hauteur)--++(0,0.5*\lgFlecheVert);
\draw [->,>=latex] (\largeur,0.5*\hauteur)--++(0,-0.5*\lgFlecheVert);
\draw (\largeur+0.4,0.5*\hauteur) node{$S_W$};
\end{tikzpicture}
\end{equation}

The mapping $D$ consists in defining the PRW from the VLMC:
the random increments of the PRW are the initial letters of a VLMC.
With the notations above, $S_n=\sum _{0\leq k\leq n}X_k$ where $X_k$ is the initial letter of $U_k$.

The mapping $L$ associates with a VLMC the process of its successive different $\alpha$-lis that turns out to be a MRC when considered together with its jump times $T_n$
-- see Section~\ref{subsec:compareSM}.
Here, $J_n^V$ is the $n$-th distinct $\alpha$-lis of the successive right-infinite words $U_0,U_1,U_2,\cdots$ and $T_n$ is the length of the $n$-th run of identical letters in the sequence $X_0,X_1X_2,\dots$
The power $V$ refers to the VLMC.

The mapping $B$ associates with a PRW $\left( S_n\right) _n$ the process of its successive different bends (changes of directions).
With our notations, $J_n^W$ is the $n$-th distinct bend and $M_n$ is the value of $S$ at the precise moment when the $n$-th bend $J_n^W$ occurs -- see Section~\ref{subsec:PRWandSM}.
The power $W$ refers to the PRW.

The mapping $S_{V}$ only consists in defining a semi-Markov process from a MRC, as stated in Section~\ref{subsec:defSM}.
The mapping $S_W$ is defined in the same manner:
it maps a MAP $\left( J_n^W,M_n\right) _n$ to the semi-Markov chain of the MRC $\left( J_n^W,M_n-M_{n-1}\right) _n$, as made precise in Definition~\ref{def:MAP}.

The mapping $N$ acts on the first coordinate by reversing words: $J_n^V=\overline{J_n^W}$.
The notation $\overline w$ stands for the reversed word of $w$:
$\overline{ab}=ba$.
For the second coordinate, remark first that $M_n-M_{n-1}$ is always of the form $k\alpha$ where $k$ is a positive integer and $\alpha$ an increment vector.
The integer $T_n$ is this $k$.

Finally, the mapping $R$ is simply the reversing of words:
$Z_n^V=\overline{Z_n^W}$.

In fact, in these particular situations (double and quadruple combs), the composition $S_V\circ L$ is a bijection -- see Remark~\ref{rem:semiMarkov}.
Therefore, all these mappings are also one-to-one, showing that all these processes are essentially equivalent.

Now that our different processes are related, let us translate the parameters,  properties and assumptions that come from the VLMC world in terms of PRW, and vice-versa.

\noindent
{\bf Dimension 1}

The PRW in dimension $1$ is driven by a VLMC based on the so-called double comb, as it was defined in Example~\ref{ex:comb}.
The contexts of this tree are the $u^kd$, which have $ud$ as an $\alpha$-lis and the $d^ku$ which have $du$ as an $\alpha$-lis ($k\geq 1$ for both families of contexts).
The cascades of the contexts write
\[
\casc\left( u^kd\right)=\prod _{j=1}^{k-1}q_{u^jd}(u)
{\rm ~~and~~}
\casc\left( d^ku\right)=\prod _{j=1}^{k-1}q_{d^ju}(d)
\]
and there are two cascade series
\[
\kappa _{ud}=\sum _{k\geq 1}\casc\left( u^kd\right)
{\rm ~~and~~}
\kappa _{du}=\sum _{k\geq 1}\casc\left( d^ku\right).
\]
Theorem~\ref{th:stablefini} guarantees that, under non-nullness assumption, this VLMC admits an invariant probability measure if, and only if $\kappa _{ud}<\infty$ and $\kappa _{du}<\infty$.
Since the double comb is a very simple context tree, one can also make a direct computation that leads to the following result:
{\it a non-null double-comb VLMC admits a $\sigma$-finite stationary measure if, and only if $\casc\left( u^kd\right)\to 0$ and $\casc\left( d^ku\right)\to 0$ when $k$ tends to infinity}.

It turns out that, on the side of the $1$-dimensional PRW, Assumption~\ref{ass:a1} as well as the expectations of the persistence times $\tau _1^u$ and $\tau _1^d$ are functions of these cascades so that one can relate the above properties of the VLMC to the results of Section~\ref{subsec:PRWdim1} on $1$-dimensional PRW.
The expectations of the waiting times are exactly the sums of cascades:
$\kappa _{ud}=\Theta _u$ and $\kappa _{du}=\Theta _d$.

Finally, one can assert:
\newcommand{\blocAss}{\hbox{Assumption \ref{ass:a1}}}
\newcommand{\blocTau}{\left(\begin{array}{c}\tau _1^u{\rm ~and~}\tau _1^d\\{\rm are~a.s.~finite}\end{array}\right)}
\newcommand{\blocCasc}{\left(\begin{array}{c}\casc\left( u^kd\right)\underset{\scriptscriptstyle k\to\infty}{\longrightarrow}0\\{\rm ~~and~~}\\[3pt]\casc\left( d^ku\right)\underset{\scriptscriptstyle k\to\infty}{\longrightarrow}0\end{array}\right)}
\newcommand{\blocMesInv}{\left(\begin{array}{c}{\rm ~The~VLMC~admits}\\{\rm a~}\sigma-{\rm finite}\\{\rm invariant~measure}\end{array}\right)}
\[
\begin{array}{ccc}
\blocCasc
&\Longleftrightarrow
&\blocAss\\[15pt]
\Updownarrow&&\Updownarrow\\[15pt]
\blocMesInv
&\Longleftrightarrow
&\blocTau
\end{array}
\]
and
\newcommand{\blocTauInt}{\left(\begin{array}{c}\tau _1^u{\rm ~and~}\tau _1^d\\{\rm are~integrable}\end{array}\right)}
\newcommand{\blocSeriesCasc}{\left(\begin{array}{c}\displaystyle\sum _{k\geq 1}\casc\left( u^kd\right)<\infty\\{\rm ~~and~~}\\[3pt]\displaystyle\sum _{k\geq 1}\casc\left( d^ku\right)<\infty\end{array}\right)}
\newcommand{\blocProbaInv}{\left(\begin{array}{c}{\rm ~The~VLMC~admits}\\{\rm a~unique~probability}\\{\rm invariant~measure}\end{array}\right)}
\[
\begin{array}{ccc}
\blocSeriesCasc
&\Longleftrightarrow
&\blocTauInt\\[15pt]
\Updownarrow&&\\[15pt]
\blocProbaInv&&
\end{array}
\]

The link between recurrence or transience of the PRW and the behaviour of the VLMC is only partial. 
For instance, the PRW may be recurrent while there is no invariant probability measure  for the VLMC.
The PRW may even be transient while the VLMC admits an invariant probability measure -- see Table \ref{tableau-rec-trans-2d}.

\noindent
{\bf Dimension 2}

The PRW in dimension $2$ is driven by a VLMC based on the so-called quadruple comb, as it is defined in Example~\ref{ex:comb}.
Here, the contexts are the $\alpha^k\beta$, where $\alpha,\beta\in\rond A=\set{\ttN,\ttE,\ttW,\ttS}$, $\alpha\neq\beta$,  $k\geq 1$.
The $\alpha$-lis of the context $\alpha^k\beta$ is $\alpha\beta$, and its cascade writes
\[
\casc\left( \alpha^k\beta\right)=\prod_{i=1}^{k-1}q_{\alpha^i\beta}(\alpha).
\]
Therefore, there are twelve cascade series, namely
\begin{equation}
\label{cascDim2}
\kappa_{\alpha \beta}=\sum_{k=1}^{\infty}
\casc\left( \alpha^k\beta\right),
~\alpha,\beta\in\rond A,~\alpha\neq\beta.
\end{equation}

As in dimension $1$, since the quadruple comb is a stable context tree having a finite set of context $\alpha$-lis, the non-null VLMC that drives the $2$-dimensional PRW admits a unique stationary probability measure if, and only if the twelve cascade series~\ref{cascDim2} converge.
This is a consequence of Theorem~\ref{th:stablefini} and, here again, due to the simplicity of the quadruple comb, on can directly check that {\it a non-null quadruple-comb VLMC admits a $\sigma$-finite stationary measure if, and only if $\casc\left( \alpha^k\beta\right)\to 0$ when $k$ tends to infinity, for every $\alpha,\beta\in\rond A$, $\alpha\neq\beta$}.

The transition matrix of the Markov process $\left( J_n\right) _n$ of the PRW bends, denoted by~$P$ in Formula~\eqref{markovsymb1}, writes also
\[
P(\beta\alpha ,\alpha\gamma)
=\sum _{n\geq 1}\casc\left( \gamma\alpha^n\beta\right)
\]
-- all other entries vanish.
Relating this expression to the definition~\eqref{defQ} of the $Q$-matrix of the VLMC leads to the following:
\begin{equation}
\label{PQ}
P(\beta\alpha ,\alpha\gamma)=Q_{\alpha\beta,\gamma\alpha}
\end{equation}
so that, up to the re-ordering that consists in reversing the indices $\alpha\beta\leadsto\beta\alpha$, the stochastic matrices $P$ and $Q$ are the same ones.
Note that, since the quadruple comb is stable, the process of the $\alpha$-lis of the VLMC is Markovian and $Q$ is its transition matrix.
Referring to the commutative diagram~\eqref{diagrammeComm}, Formula~\eqref{PQ} amounts to saying that the Markov chains $\left( J_n^V\right) _n$ and $\left(\overline{J_n^W}\right) _n$ are identical.

In terms of persistence times of the PRW vs stationary measures for the VLMC,  the properties stated in Section~\ref{subsec:PRWdim2} show that the following equivalences hold.
\newcommand{\blocAssDeux}{\hbox{Assumption \ref{ass:a2}}}
\newcommand{\blocTn}{\left(\begin{array}{c}\forall n,~T_n{\rm ~is~a.s.~finite}\end{array}\right)}
\newcommand{\blocCascDeux}{\left(\begin{array}{c}{\rm for~all~}\alpha,\beta\in\rond A,\alpha\neq\beta ,\\[8pt]\casc\left( \alpha^k\beta\right)\underset{\scriptscriptstyle k\to\infty}{\longrightarrow}0\end{array}\right)}
\newcommand{\blocMesInvDeux}{\left(\begin{array}{c}{\rm ~The~VLMC~admits}\\{\rm a~}\sigma-{\rm finite}\\{\rm invariant~measure}\end{array}\right)}
\[
\begin{array}{ccc}
\blocCascDeux
&\Longleftrightarrow
&\blocAssDeux\\[15pt]
\Updownarrow&&\Updownarrow\\[15pt]
\blocMesInvDeux
&\Longleftrightarrow
&\blocTn
\end{array}
\]
and
\newcommand{\blocTnInt}{\left(\forall n,~T_n{\rm ~is~integrable}\right)}
\newcommand{\blocSeriesCascDeux}{\left(\begin{array}{c}{\rm for~all~}\alpha,\beta\in\rond A,\alpha\neq\beta ,\\[8pt]\displaystyle\sum _{k\geq 1}\casc\left( \alpha^k\beta\right)<\infty\end{array}\right)}
\newcommand{\blocProbaInvDeux}{\left(\begin{array}{c}{\rm ~The~VLMC~admits}\\{\rm a~unique~probability}\\{\rm invariant~measure}\end{array}\right)}
\[
\begin{array}{ccc}
\blocSeriesCascDeux
&\Longleftrightarrow
&\blocTnInt\\[15pt]
\Updownarrow&&\\[15pt]
\blocProbaInvDeux&&
\end{array}
\]
The counterexample cited in Theorem~\ref{conjecture} is enlighted by these equivalences:
an example of recurrent $2$-dimensional PRW having a transient  skeleton $\left( M_n\right) _n$ cannot be found without assuming that the $T_n$ are a.s. finite but non-integrable, as shown in~\cite{cenac:hal-01658494}.
Reading the above equivalences shows that such a PRW must be driven by a VLMC the series of cascades of which diverge while their general terms tend to zero at infinity.

\bibliographystyle{agsm} 
\bibliography{paccaut}

\vfill\pagebreak
\
\thispagestyle{empty}
\end{document}